\documentclass[11pt]{article}

\usepackage[utf8]{inputenc}
\usepackage[T1]{fontenc}
\usepackage[english]{babel}
\usepackage{lmodern}
\usepackage{amsmath,amssymb,amsthm,mathtools}
\usepackage[margin=2.6cm]{geometry}
\usepackage{enumitem}
\usepackage[numbers,sort&compress]{natbib}
\usepackage{xcolor}
\usepackage[colorlinks=true,linkcolor=blue!50!black,
            citecolor=green!40!black,urlcolor=blue!50!black]{hyperref}

\numberwithin{equation}{section}

\newtheorem{theorem}{Theorem}[section]
\newtheorem{proposition}[theorem]{Proposition}
\newtheorem{lemma}[theorem]{Lemma}
\newtheorem{corollary}[theorem]{Corollary}
\theoremstyle{definition}
\newtheorem{definition}[theorem]{Definition}
\newtheorem{assumption}[theorem]{Assumption}
\newtheorem*{hypH}{Assumption (H)}
\newtheorem*{hypHg}{Assumption (H$_\gamma$)}
\newtheorem{example}[theorem]{Example}
\theoremstyle{remark}
\newtheorem{remark}[theorem]{Remark}

\newcommand{\R}{\mathbb{R}}

\newcommand{\E}{\mathbb{E}}
\newcommand{\Prob}{\mathbb{P}}
\newcommand{\iu}{\mathrm{i}}
\newcommand{\dd}{\mathrm{d}}
\newcommand{\Law}{\mathcal{L}}
\newcommand{\cP}{\mathcal{P}}
\newcommand{\cA}{\mathcal{A}}
\newcommand{\cF}{\mathcal{F}}
\newcommand{\Sch}{\mathcal{S}}
\newcommand{\Wr}{W_\rho}
\newcommand{\one}{\mathbf{1}}

\newcommand{\norm}[1]{\left\|#1\right\|}
\newcommand{\ip}[2]{\langle #1,\, #2\rangle}

\title{McKean--Vlasov flights without moments: well-posedness,\\
propagation of chaos and a critical threshold for law-dependent jumps}
\author{Houssine El Jeddaoui\\[0.8em]
\small Laboratoire de Math\'ematiques Informatique et Applications (LAMIA),\\
\small Universit\'e des Antilles, Pointe-\`a-Pitre, Guadeloupe, France.\\[0.4em]
\small E-mail:
\href{mailto:houssine.eljeddaoui@univ-antilles.fr}{houssine.eljeddaoui@univ-antilles.fr}}
\date{}

\begin{document}
\maketitle

\begin{abstract}
We study McKean--Vlasov equations driven by an isotropic $2s$-stable L\'evy
noise, $s\in(0,1)$, whose solutions may have no finite mean, with
coefficients that are Lipschitz with respect to a bounded Wasserstein-type
distance. When the jump noise is additive, we prove strong well-posedness
for every initial law and every $s$, without any moment assumption, and
propagation of chaos at the rate $N^{-1/2}$ for kernel interactions. When
the size of the jumps depends on the law, we prove the same results for
$s>\frac12$, and we show by an example that uniqueness fails for every
$s<\frac12$. The value $s=\frac12$ thus separates well-posedness from
non-uniqueness; the critical case itself remains open. These results
answer, for stable noise, the question of uniqueness without moments raised
by Jourdain, M\'el\'eard and Woyczy\'nski and by Cavallazzi.
\end{abstract}

\textbf{Keywords.} McKean--Vlasov equations; stable processes;
L\'evy flights; propagation of chaos; nonlinear Fokker--Planck equations;
fractional Laplacian; pseudo-differential operators.

\textbf{MSC 2020.} 60H10, 60K35, 60G52, 60J76, 35Q84, 35R11, 47G30.

\section{Introduction}

Stable processes, or L\'evy flights, are the standard model for random
displacements with heavy tails. They describe anomalous diffusion in physics
\cite{MetzlerKlafter2000}, foraging patterns in ecology \cite{Viswanathan1996}
and price fluctuations in finance \cite{Mandelbrot1963}. When many such
agents interact through their empirical distribution, the mean-field limit
is a McKean--Vlasov equation, whose coefficients depend on the law of the
solution itself \cite{McKean1966,Sznitman1991}. We call its solutions
\emph{McKean--Vlasov flights}.
In this paper we consider
\begin{equation}\label{eq:intro-MV0}
  \dd X_t=b(X_t,\mu_t)\,\dd t+\sqrt{2\sigma(X_t,\mu_t)}\,\dd W_t
  +\gamma(X_{t-},\mu_t)\,\dd L_t,\qquad \mu_t=\Law(X_t),
\end{equation}
where $W$ is a Brownian motion and $L$ an isotropic $2s$-stable process,
$s\in(0,1)$. The coefficient $\gamma$ is the size of the jumps; the case
$\gamma\equiv I_d$ is that of additive jump noise. The law $\mu_t$ solves a
nonlinear Fokker--Planck equation whose generator is a nonlocal operator of
fractional order $2s$.

Letting the size of the jumps depend on the population is natural in
models of interacting agents. A forager may make longer flights where its
congeners are sparse, and the size of price jumps depends on the state of
the market. In such models the displacements are heavy-tailed, the
population may have no mean position, and the law enters the noise itself.
This is where the existing theory, which relies on moments, does not
apply, and this is the situation we study.

For Brownian noise, the basic theory is classical and requires no moments.
McKean's theorem, in the form proved by Sznitman \cite{Sznitman1991}, gives
existence, uniqueness and propagation of chaos at rate $N^{-1/2}$ for every
initial law \cite[Thm.~3.1]{ChaintronDiez2022}. For stable noise, moments
are a real issue: a $2s$-stable law has no moment of order $2s$, and no mean
when $s\le\frac12$. Yet all existing uniqueness results for
\eqref{eq:intro-MV0} restrict the initial law, except when the jumps are
additive and the law enters through a convolution, as in the work of Hao
\cite{Hao2025} and of Hao, Ren and Wu \cite{HaoRenWu2024}. They are formulated in Wasserstein spaces
with moments of order at least one
\cite{JMW2008,Cavallazzi2025,BLW2024,NeelimaEtAl,DengHuang2023}, which excludes
$s\le\frac12$, or with a moment of small order and a non-degenerate noise
\cite{DengHuang2025}, or they require an initial density and $s>\frac12$
\cite{BarbuRockner2024}. Whether uniqueness holds without moments is stated
as an open question by Jourdain, M\'el\'eard and Woyczy\'nski
\cite[Rem.~8]{JMW2008} and by Cavallazzi \cite{Cavallazzi2025}. The
question we address is therefore the following: \emph{under which
conditions does the McKean--Sznitman theory hold for McKean--Vlasov flights,
without any moment assumption?}

Our answer is as follows. We assume that $b$, $\sqrt{2\sigma}$ and $\gamma$ are Lipschitz in the
position and in the law, the law being measured by the bounded
Wasserstein-type distance $\Wr$, which is defined on all probability
measures. Under this assumption alone, we prove the following.
\begin{enumerate}[label=(\roman*)]
  \item If the jump noise is additive, the theory holds for every
    $s\in(0,1)$: strong existence and pathwise uniqueness for every initial
    law, the nonlinear Fokker--Planck equation, and propagation of chaos at
    the rate $N^{-1/2}$ for kernel interactions, with a constant independent
    of the initial law and of $s$ (Theorems~\ref{thm:main-wp}
    and~\ref{thm:main-chaos}).
  \item If the size of the jumps depends on the law, the theory holds for
    every $s\in(\frac12,1)$, again for every initial law
    (Theorem~\ref{thm:main-jump}).
  \item If the size of the jumps depends on the law and $s<\frac12$,
    uniqueness fails in general: we give a bounded coefficient $\gamma(\mu)$,
    Lipschitz for $\Wr$, for which the equation has two solutions with
    different laws (Proposition~\ref{prop:main-counter}).
\end{enumerate}
In short, the paper shows when the McKean--Sznitman
theory survives heavy tails without moments: always for additive jumps, for
law-dependent jumps above the threshold $s=\frac12$, and in general not
below it, where the jumps have no mean.

To the best of our knowledge, Theorem~\ref{thm:main-jump} is the first
well-posedness result for McKean--Vlasov equations with law-dependent jumps
of infinite activity that requires neither moments nor densities on the
initial law. Cavallazzi \cite{Cavallazzi2025}, Bao, Liu and Wang
\cite{BLW2024} and Neelima et al.\ \cite{NeelimaEtAl} require moments of order
at least one, Deng and Huang \cite{DengHuang2023} moments of order at least
one, a non-degenerate noise and a jump coefficient independent of the
position, Deng and Huang \cite{DengHuang2025} a moment of small order, a
non-degenerate noise and $s\le\frac12$, Barbu and R\"ockner \cite{BarbuRockner2024} an
initial density, and Graham \cite{Graham1992b} a bounded jump rate. Proposition~\ref{prop:main-counter}
shows that the restriction $s>\frac12$ cannot be relaxed to any
$s<\frac12$; only the critical case $s=\frac12$ remains open. Unlike the counterexamples of
Cavallazzi \cite{Cavallazzi2025} and of Deng and Huang \cite{DengHuang2025},
it uses a bounded coefficient that is
Lipschitz for the bounded distance $\Wr$, so that non-uniqueness is caused
by the heavy tails of the noise and not by an unbounded dependence on the
law. For additive jumps, well-posedness without moments is known when the
law enters through a convolution, even for H\"older or singular kernels
\cite{Hao2025,HaoRenWu2024}. Our result (i) requires Lipschitz regularity,
but allows a general dependence on the law and a law-dependent, possibly
degenerate Brownian part, and it gives the rate $N^{-1/2}$ with a constant
independent of the initial law, which seems to be new for stable noise. As
a consequence, our results answer the uniqueness question of
\cite{JMW2008,Cavallazzi2025} for stable noise, for additive jumps and for
law-dependent jumps with $s>\frac12$, and show that it has a negative answer
in general when $s<\frac12$.

We believe these results are of interest for several reasons. First, they show that moment conditions are not intrinsic to the
McKean--Vlasov theory with jumps. They come from the use of Wasserstein
distances and disappear with a bounded distance, which is the natural
choice for heavy-tailed laws. Indeed, even if $\mu_0$ has moments, the law
$\mu_t$ has in general no moment of order $2s$ (for instance when the jumps
are additive and $b$, $\sigma$ are bounded), so the natural state space of
McKean's nonlinear semigroup $\mu_0\mapsto\mu_t$ is the set of all
probability measures. In particular, mean-field limits of L\'evy
flights with infinite mean are well posed and are approximated by particle
systems at the same rate as in the Brownian case. Second, they show where heavy tails
become an obstruction: never when the jumps are additive and, when their
size depends on the law, below the threshold $s=\frac12$ at which the
jumps lose their mean. Third, the proofs are short and elementary and do not
use regularisation by the noise, so they apply to degenerate equations,
such as the pure flight $\dd X_t=b(X_t,\mu_t)\,\dd t+\dd L_t$. Finally, we describe the generator through its nonlinear symbol, in the
spirit of the theory of pseudo-differential operators
\cite{Courrege1966,Jacob2001,BSW2013,SchillingSchnurr2010}.

\subsection{Main results}

Let $s\in(0,1)$, let $W$ be a standard Brownian motion in $\R^d$ and let $L$
be an independent isotropic $2s$-stable L\'evy process, with
$\E[e^{\iu\xi\cdot L_t}]=e^{-t|\xi|^{2s}}$. We consider
\begin{equation}\label{eq:intro-MV}
  \dd X_t=b(X_t,\mu_t)\,\dd t+\Sigma(X_t,\mu_t)\,\dd W_t+\dd L_t,
  \qquad \mu_t=\Law(X_t),\qquad \Sigma:=\sqrt{2\sigma}.
\end{equation}
On the set $\cP(\R^d)$ of all probability measures on $\R^d$ we use the
distance
\[
  \Wr(\mu,\nu):=\inf_{\eta\in\Pi(\mu,\nu)}\int\big(|x-y|\wedge1\big)\,
  \eta(\dd x,\dd y),
\]
where $\Pi(\mu,\nu)$ is the set of couplings of $\mu$ and $\nu$. This distance
metrises weak convergence and is finite without any moment
condition. Our standing assumption is the following.

\begin{hypH}
There is $K<\infty$ such that, for all $x,y\in\R^d$ and $\mu,\nu\in\cP(\R^d)$,
\[
  |b(x,\mu)-b(y,\nu)|+|\Sigma(x,\mu)-\Sigma(y,\nu)|
  \le K\big(|x-y|+\Wr(\mu,\nu)\big).
\]
\end{hypH}

For $\mu\in\cP(\R^d)$ we write
$\cA_\mu f:=b(\cdot,\mu)\cdot\nabla f+\sigma(\cdot,\mu)\Delta f-(-\Delta)^sf$.

\begin{theorem}\label{thm:main-wp}
Let $s\in(0,1)$ and assume \textup{(H)}. For every $\mu_0\in\cP(\R^d)$ and
every $T>0$, the following hold.
\begin{enumerate}[label=(\roman*)]
  \item Equation \eqref{eq:intro-MV} with $\Law(X_0)=\mu_0$ has a strong
    solution on $[0,T]$, and any two solutions are indistinguishable. The
    flow $t\mapsto\mu_t$ is $\Wr$-continuous and depends only on $\mu_0$.
  \item For every $\varphi\in C_c^\infty(\R^d)$, the process
    $\displaystyle \varphi(X_t)-\varphi(X_0)-\int_0^t\cA_{\mu_r}\varphi(X_r)\,\dd r$ is a
    martingale, and $(\mu_t)$ is a weak solution of the nonlinear
    Fokker--Planck equation $\partial_t\mu_t=\cA_{\mu_t}^*\mu_t$, that is,
    \[
      \ip{\mu_t}{\varphi}=\ip{\mu_0}{\varphi}+\int_0^t\ip{\mu_r}{\cA_{\mu_r}\varphi}\,\dd r .
    \]
  \item The maps $T_t(\mu_0):=\mu_t$ form a semigroup of nonlinear maps on
    $\cP(\R^d)$.
\end{enumerate}
\end{theorem}

Let $(X^i_0,W^i,L^i)_{i\ge1}$ be i.i.d.\ copies of $(X_0,W,L)$, and consider
the particle system
\[
  \dd X^{i,N}_t=b(X^{i,N}_t,\mu^N_t)\,\dd t+\Sigma(X^{i,N}_t,\mu^N_t)\,\dd W^i_t+\dd L^i_t,
  \qquad \mu^N_t:=\frac1N\sum_{j=1}^N\delta_{X^{j,N}_t},
\]
with $X^{i,N}_0=X^i_0$. Let $\bar X^i$ be the solution of
\eqref{eq:intro-MV} driven by $(X^i_0,W^i,L^i)$.

\begin{theorem}\label{thm:main-chaos}
Let $s\in(0,1)$, assume \textup{(H)} and let $\mu_0\in\cP(\R^d)$. Then, for
every $T>0$ and $i\ge1$,
\[
  \lim_{N\to\infty}\E\Big[\sup_{t\le T}|X^{i,N}_t-\bar X^i_t|^2\Big]=0 .
\]
If moreover $\displaystyle b(x,\mu)=\int\beta(x,y)\,\mu(\dd y)$ and
$\displaystyle \Sigma(x,\mu)=\int\varsigma(x,y)\,\mu(\dd y)$, where
$\beta:\R^d\times\R^d\to\R^d$ and $\varsigma:\R^d\times\R^d\to[0,\infty)$
are bounded and Lipschitz, then
\[
  \E\Big[\sup_{t\le T}|X^{i,N}_t-\bar X^i_t|^2\Big]\le\frac CN,
\]
where $C$ depends only on $d$, $T$ and the kernels, but not on $\mu_0$ or $s$.
\end{theorem}

Theorems~\ref{thm:main-wp} and~\ref{thm:main-chaos} are proved in
Section~\ref{sec:nonlinear} (Theorems~\ref{thm:wellposed}, \ref{thm:chaos},
\ref{thm:chaos-rate}, \ref{thm:MP} and
Proposition~\ref{prop:nl-semigroup}).

We now let the jumps depend on the law, that is, we consider
\eqref{eq:intro-MV0} with a general coefficient
$\gamma:\R^d\times\cP(\R^d)\to\R^{d\times d}$, under the assumption
\begin{hypHg}
Assumption (H) holds, and $|\gamma(x,\mu)-\gamma(y,\nu)|\le K\big(|x-y|+\Wr(\mu,\nu)\big)$ for all $x,y\in\R^d$ and $\mu,\nu\in\cP(\R^d)$.
\end{hypHg}
The generator is then the nonlocal operator \eqref{eq:gen-jump}, with
symbol
$p^\gamma(x,\xi;\mu)=-\iu b(x,\mu)\cdot\xi+\sigma(x,\mu)|\xi|^2+|\gamma(x,\mu)^\top\xi|^{2s}$,
in which the law enters the jump part.

\begin{theorem}\label{thm:main-jump}
Let $s\in(\frac12,1)$, assume \textup{(H$_\gamma$)} and let
$\mu_0\in\cP(\R^d)$. Then \eqref{eq:intro-MV0} has a unique strong solution
on $[0,T]$ for every $T>0$, the map $\mu_0\mapsto\mu_t$ is a nonlinear Markov
semigroup, and, if the coefficients are bounded, $(\mu_t)$ solves the
nonlinear Fokker--Planck equation. Moreover, for the particle
system \eqref{eq:particles} and every $p\in[1,2s)$,
$\lim_{N\to\infty}\E[\sup_{t\le T}|X^{i,N}_t-\bar X^i_t|^p]=0$, and for
interactions through bounded Lipschitz kernels,
\[
  \Big(\E\Big[\sup_{t\le T}|X^{i,N}_t-\bar X^i_t|^p\Big]\Big)^{1/p}
  \le\frac{C}{\sqrt N},
\]
with $C$ independent of $\mu_0$.
\end{theorem}

\begin{proposition}\label{prop:main-counter}
Let $s\in(0,\frac12)$ and $\displaystyle \gamma(\mu):=\Bigl(\int(|x|\wedge1)\,\mu(\dd x)\Bigr)I_d$, which
satisfies $|\gamma(\mu)-\gamma(\nu)|\le\sqrt d\,\Wr(\mu,\nu)$. The equation
$\dd X_t=\gamma(\mu_t)\,\dd L_t$, $X_0=0$, has two solutions with
different laws.
\end{proposition}

Theorem~\ref{thm:main-jump} and Proposition~\ref{prop:main-counter} are
proved in Section~\ref{sec:jumps} (Theorems~\ref{thm:jump-wp}
and~\ref{thm:jump-chaos}, Proposition~\ref{prop:counter}).

Three features of the assumptions should be noted. First, in
Theorems~\ref{thm:main-wp} and~\ref{thm:main-chaos} the law enters only
through the drift and the Brownian coefficient; by
Proposition~\ref{prop:main-counter}, this restriction is necessary for
$s<\frac12$. Second, the distance $\Wr$ is defined on all of
$\cP(\R^d)$ and $(\cP(\R^d),\Wr)$ is complete. Assumption (H) holds for
interactions through bounded Lipschitz kernels (Example~\ref{ex:kernel}),
and excludes interactions such as $\displaystyle \int(x-y)\,\mu(\dd y)$, which are not
defined when the mean is infinite. Third, no non-degeneracy is required:
$\sigma\equiv0$ is allowed, which gives the pure McKean--Vlasov flight
$\dd X_t=b(X_t,\mu_t)\,\dd t+\dd L_t$.

We now describe the main ideas of the proofs. They follow the synchronous coupling of McKean and Sznitman. The
point specific to L\'evy noise is that two solutions driven by the same
$(X_0,W,L)$ have a continuous difference, since the additive jumps cancel.
$L^2$ estimates therefore apply to this difference, although the solutions
themselves have no moment of order $2s$ (Lemma~\ref{lem:stability}).
Together with the completeness of $(\cP(\R^d),\Wr)$, this yields a Picard
fixed point for every initial law. For the rate in
Theorem~\ref{thm:main-chaos}, the fluctuations of the empirical measure are
controlled by conditional independence, with constants that do not involve
the noise.

For law-dependent jumps, the difference of two solutions has jumps, and the
$L^2$ argument fails. When $s>\frac12$, we work instead in $L^p$ with
$p\in(1,2s)$: small jumps are controlled through their predictable quadratic
variation and Lenglart's inequality, and large jumps in $L^p$ through their compensator
(Lemma~\ref{lem:Lp}). Since $p\ge1$, the bound
$\Wr\le(\E|D|^p)^{1/p}$ is linear, which gives a contraction on a short
time interval whose length does not depend on the initial law. The
solution is then obtained by concatenation. When $s<\frac12$, only moments
of order $p<1$ are available and the bound becomes
$\Wr\lesssim\bar w^{\,p}$, which is not a contraction. The counterexample
exploits exactly this: a deterministic time change $S_t$ of the stable
process solves an ordinary differential equation $\dot S=g(S)$ with
$g(S)\ge cS^{2s}$ near $0$, which has two solutions starting from $0$.

We end this introduction with a more detailed account of the literature. McKean--Vlasov equations with jumps go back to Graham
\cite{Graham1992,Graham1992b}, who treats jumps with a bounded rate; see
also M\'el\'eard \cite{Meleard1996} and the review \cite{ChaintronDiez2022}.
For L\'evy noise of infinite activity, Jourdain, M\'el\'eard and
Woyczy\'nski \cite{JMW2008} consider $\dd X_t=\sigma(X_{t-},\mu_t)\,\dd Z_t$,
in connection with fractal Burgers equations \cite{BFW1998}. They prove
well-posedness and propagation of chaos for square-integrable $X_0$ and $Z$,
and existence without integrability. Cavallazzi \cite{Cavallazzi2025}
proves well-posedness and quantitative chaos with moments of order
$\beta\in[1,2]$, Bao, Liu and Wang \cite{BLW2024} replace the Lipschitz
condition by monotonicity, Neelima et al.\ \cite{NeelimaEtAl} work in
$\cP_2$, Erny \cite{Erny2022} treats locally Lipschitz coefficients with
exponential moments, and the long-time behaviour is studied in
\cite{LMW2021}. For stable noise, regularisation by the noise allows H\"older
or singular drifts: Wu and Hao \cite{WuHao2023} and Frikha, Konakov and
Menozzi \cite{FKM2021} treat $s>\frac12$ and $s>\frac13$ respectively,
Chaudru de Raynal, Jabir and Menozzi \cite{CdRJM2025} and Hao, R\"ockner and
Zhang \cite{HaoRocknerZhang} distributional or singular kernels, and Hao
\cite{Hao2025} and Hao, Ren and Wu \cite{HaoRenWu2024} convolution drifts for
every $s$ and every initial law. Deng and Huang treat non-degenerate jump coefficients depending on the
law. In \cite{DengHuang2023}, $s\in(\frac12,1)$, the coefficient
$\sigma_t(\mu)$ does not depend on the position, and the initial law is in
$\cP_k$ with $k\in[1,2s)$. In \cite{DengHuang2025}, $s\le\frac12$, the
coefficient depends on the position and the law, the initial law is in
$\cP_\eta$ with $\eta\in(0,2s)$, and the drift may depend on the law only
when $s>\frac14$.
Barbu and R\"ockner \cite{BarbuRockner2024} treat coefficients depending
pointwise on the density, for $s\in(\frac12,1)$ and initial densities in
$L^1\cap L^2$, with uniqueness in law among solutions with bounded densities.

The paper is organised as follows. Section~\ref{sec:prelim} collects the notation, the setting and
the standing assumptions. Section~\ref{sec:nonlinear} treats additive jump noise.
Besides the proofs of Theorems~\ref{thm:main-wp} and~\ref{thm:main-chaos},
it contains the carr\'e du champ of $\cA_\mu$ (Theorem~\ref{thm:MP}), the
nonlinear symbol $p(x,\xi;\mu)$ with a probabilistic formula
(Theorem~\ref{thm:nl-symbol}), and the case of coefficients depending on
$\mu$ only (Corollary~\ref{cor:explicit}), with an example of a mean-field
flight whose centre is driven to the origin although its mean does not
exist (Example~\ref{ex:flight}). Section~\ref{sec:jumps} treats
law-dependent jumps and contains the counterexample.
Appendix~\ref{sec:linear} recalls, with a single sign convention, the
correspondence between generator, symbol and law for the linear operator
$b\cdot\nabla+\sigma\Delta-(-\Delta)^s$
(Propositions~\ref{prop:ito}--\ref{prop:schilling}).

\section{Notation, setting and standing assumptions}\label{sec:prelim}

Throughout, $d \ge 1$, $e_\xi(x) := e^{\iu\xi\cdot x}$ for $\xi,x \in \R^d$, and $\cP(\R^d)$ is the set of Borel probability measures on $\R^d$. The law of a random variable $X$ is denoted by $\Law(X)$. For $\mu \in \cP(\R^d)$ and a bounded measurable $\varphi$, we write
\[
\ip{\mu}{\varphi} := \int \varphi \, \dd\mu.
\]

For $f \in L^1(\R^d)$, the Fourier transform is $\displaystyle \hat{f}(\xi) := \int_{\R^d} f(x) e^{-\iu\xi\cdot x} \, \dd x$.

On the Schwartz space $\Sch(\R^d)$, the inverse and the derivative property are given by $\displaystyle f(x) = (2\pi)^{-d} \int_{\R^d} \hat{f}(\xi) e^{\iu\xi\cdot x} \, \dd\xi$ and $\widehat{\partial_j f}(\xi) = \iu\xi_j \hat{f}(\xi)$.

For a probability measure we use the \emph{characteristic function}
\[
\chi_\mu(\xi) := \ip{\mu}{e_\xi} = \int_{\R^d} e^{\iu\xi\cdot x} \, \mu(\dd x),
\]
so that, for $\mu(\dd x) = f(x) \, \dd x$, $\chi_\mu(\xi) = \hat{f}(-\xi)$.

For two nonnegative quantities $f$ and $g$, we write $f\lesssim g$ if $f\le Cg$ for a constant $C>0$ that does not depend on the variables involved.

All L\'evy processes are taken with respect to a filtration $(\cF_t)$ satisfying the usual conditions: $L_0=0$, $L$ has c\`adl\`ag paths, and $L_t-L_s$ is independent of $\cF_s$ and distributed as $L_{t-s}$ for $s<t$.

\begin{theorem}[L\'evy--Khintchine {\cite[Thms.~7.10, 8.1]{Sato1999}}]
\label{thm:LK}
Let $L$ be a L\'evy process. There exists a unique triplet $(b,Q,\nu)$, with
$b\in\R^d$, $Q$ symmetric non-negative definite and $\nu$ a measure on
$\R^d\setminus\{0\}$ with $\displaystyle \int(1\wedge|y|^2)\,\nu(\dd y)<\infty$, such that
\begin{equation}\label{eq:LK}
  \E\big[e^{\iu\xi\cdot L_t}\big]=e^{-t\psi(\xi)},\qquad
  \psi(\xi)=-\iu\, b\cdot\xi+\tfrac12\,\xi\cdot Q\xi
  +\int_{\R^d\setminus\{0\}}\big(1-e^{\iu\xi\cdot y}
  +\iu\,\xi\cdot y\,\one_{\{|y|<1\}}\big)\,\nu(\dd y).
\end{equation}
Conversely, every such triplet arises from a L\'evy process. The function
$\psi$ is called the \emph{characteristic exponent} of $L$.
\end{theorem}

\begin{example}\label{ex:basic}
\begin{enumerate}[label=(\alph*)]
  \item \emph{Drift.} $L_t=bt$ has triplet $(b,0,0)$ and
    $\psi(\xi)=-\iu\, b\cdot\xi$.
  \item \emph{Brownian motion.} If $W$ is a standard Brownian motion in
    $\R^d$ and $\sigma\ge0$, then $L_t=\sqrt{2\sigma}\,W_t$ has triplet
    $(0,2\sigma I_d,0)$ and $\psi(\xi)=\sigma|\xi|^2$.
  \item \emph{Isotropic stable process.} See Definition~\ref{def:stable}.
\end{enumerate}
\end{example}

Fix $s\in(0,1)$ and set
\begin{equation}\label{eq:cds}
  c_{d,s}:=\frac{4^s\,\Gamma(\tfrac d2+s)}{\pi^{d/2}\,|\Gamma(-s)|},
  \qquad
  \nu_s(\dd y):=\frac{c_{d,s}}{|y|^{d+2s}}\,\dd y .
\end{equation}
Then $\nu_s$ is a symmetric L\'evy measure, and the constant is chosen so that
\begin{equation}\label{eq:cos}
  \int_{\R^d\setminus\{0\}}\big(1-\cos(\xi\cdot y)\big)\,\nu_s(\dd y)
  =|\xi|^{2s},\qquad \xi\in\R^d,
\end{equation}
see \cite[Sec.~3]{DNPV2012} and \cite{Kwasnicki2017}.

\begin{definition}\label{def:fraclap}
For $f\in C_b^2(\R^d)$ (bounded with bounded derivatives up to order two),
\begin{equation}\label{eq:fraclap}
  -(-\Delta)^s f(x):=\int_{\R^d\setminus\{0\}}
  \big[f(x+y)-f(x)-\nabla f(x)\cdot y\,\one_{\{|y|<1\}}\big]\,\nu_s(\dd y).
\end{equation}
The integrand is bounded by
$\min\big(\tfrac12\norm{D^2f}_\infty|y|^2,\,2\norm{f}_\infty+\norm{\nabla f}_\infty\big)$,
so $(-\Delta)^sf$ is bounded and continuous. By symmetry of $\nu_s$, the
compensating term may be replaced by a principal value.
\end{definition}

Applying \eqref{eq:fraclap} to $e_\xi$ and using \eqref{eq:cos} and the
symmetry of $\nu_s$ gives
\begin{equation}\label{eq:fraclap-plane}
  -(-\Delta)^s e_\xi=-|\xi|^{2s}\,e_\xi .
\end{equation}
As a consequence, for $f\in\Sch(\R^d)$,
$\widehat{(-\Delta)^sf}(\xi)=|\xi|^{2s}\hat f(\xi)$, which recovers the
spectral definition \cite{DNPV2012,Kwasnicki2017}.

\begin{definition}\label{def:stable}
The \emph{isotropic $2s$-stable L\'evy process} $L=L^{(s)}$ is the L\'evy
process with triplet $(0,0,\nu_s)$. By \eqref{eq:LK}, \eqref{eq:cos} and
symmetry, $\psi(\xi)=|\xi|^{2s}$, that is,
$\E[e^{\iu\xi\cdot L_t}]=e^{-t|\xi|^{2s}}$. We write $\pi_t:=\Law(L_t)$.
\end{definition}

\begin{remark}\label{rem:moments}
Since $\nu_s(\{|y|>r\})=C\,r^{-2s}$, one has $\E|L_t|^q<\infty$ if and only
if $q<2s$ \cite[Thm.~25.3]{Sato1999}. In particular, $L$ has infinite
variance, and infinite mean when $s\le\tfrac12$. This is why we avoid
moment assumptions in Sections~\ref{sec:nonlinear} and~\ref{sec:jumps}. Brownian motion
corresponds formally to the limit $s=1$: since $c_{d,s}\to0$ as
$s\uparrow1$, the jump measure degenerates and only the Gaussian part
survives.
\end{remark}

By the L\'evy--It\^o decomposition, and since the triplet is $(0,0,\nu_s)$,
\begin{equation}\label{eq:LI}
  L_t=\int_0^t\!\!\int_{\{|y|<1\}}y\,\tilde N(\dd r,\dd y)
  +\int_0^t\!\!\int_{\{|y|\ge1\}}y\,N(\dd r,\dd y),
\end{equation}
where $N$ is the Poisson random measure of the jumps of $L$, with intensity
$\dd r\,\nu_s(\dd y)$, and $\tilde N:=N-\dd r\,\nu_s(\dd y)$.

By Courr\`ege's theorem \cite{Courrege1966} (see \cite[Ch.~2]{BSW2013}), a
linear operator $\cA:C_c^\infty(\R^d)\to C(\R^d)$ satisfying the positive
maximum principle can be written as
\begin{equation}\label{eq:pdo}
  \cA f(x)=-\frac{1}{(2\pi)^d}\int_{\R^d}p(x,\xi)\,\hat f(\xi)\,
  e^{\iu\xi\cdot x}\,\dd\xi,
  \qquad f\in C_c^\infty(\R^d),
\end{equation}
where, for every $x$, $\xi\mapsto p(x,\xi)$ is of the form \eqref{eq:LK} with
an $x$-dependent triplet, plus a killing term $c(x)\ge0$.

\begin{definition}\label{def:symbol}
Let $\cA$ be a linear operator defined on $C_b^2(\R^d)$, viewed as acting on
complex-valued functions by linearity. Its \emph{symbol} is
\begin{equation}\label{eq:symbol}
  p(x,\xi):=-\,e_{-\xi}(x)\,(\cA e_\xi)(x),\qquad x,\xi\in\R^d .
\end{equation}
\end{definition}

With this sign, the symbol of the generator of a L\'evy process is its
characteristic exponent, $p(x,\xi)=\psi(\xi)$
(Proposition~\ref{prop:analytic}), and \eqref{eq:pdo} reads
$\cA=-p(x,D)$, in agreement with \cite{BSW2013,SchillingSchnurr2010}.

Let $\rho(x,y):=|x-y|\wedge 1$. For $\mu,\nu\in\cP(\R^d)$, define
\begin{equation}\label{eq:Wrho}
  \Wr(\mu,\nu):=\inf_{\eta\in\Pi(\mu,\nu)}\int\rho(x,y)\,\eta(\dd x,\dd y),
\end{equation}
where $\Pi(\mu,\nu)$ is the set of couplings of $\mu$ and $\nu$.

\begin{lemma}\label{lem:Wrho}
\begin{enumerate}[label=(\alph*)]
  \item $(\cP(\R^d),\Wr)$ is a complete separable metric space, and $\Wr$
    metrises weak convergence.
  \item If $X\sim\mu$ and $Y\sim\nu$ are defined on the same probability
    space, then
    $\Wr(\mu,\nu)\le\E[|X-Y|\wedge1]\le(\E|X-Y|^2)^{1/2}$.
  \item If $f$ is bounded and Lipschitz, then
    $|\ip{\mu}{f}-\ip{\nu}{f}|\le\max(\mathrm{Lip}(f),2\norm{f}_\infty)\,\Wr(\mu,\nu)$.
  \item $\Wr\le1$.
\end{enumerate}
\end{lemma}

\begin{proof}
(a) Since $(\R^d,\rho)$ is a complete, separable and bounded metric space,
$\Wr$ is its $1$-Wasserstein distance; it is complete and metrises weak
convergence by \cite[Thms.~6.9 and~6.18]{Villani2009}.

(b) The joint law of $(X,Y)$ is a coupling, which gives the first
inequality; the second is Jensen's inequality.

(c) For all $x,y\in\R^d$,
\[
  |f(x)-f(y)|\le\min\bigl(\mathrm{Lip}(f)|x-y|,\,2\norm{f}_\infty\bigr)
  \le\max\bigl(\mathrm{Lip}(f),\,2\norm{f}_\infty\bigr)\,\rho(x,y),
\]
and we integrate against an optimal coupling.

(d) This follows from $\rho\le1$.
\end{proof}

We now fix the setting used in the rest of the paper. Matrices are given the Frobenius norm $|A|:=\bigl(\sum_{i,j}A_{ij}^2\bigr)^{1/2}$. Let $(\Omega,\cF,(\cF_t),\Prob)$ be a filtered probability space satisfying the usual conditions, carrying an $\cF_0$-measurable random variable $X_0$ with law $\mu_0$, a standard $d$-dimensional Brownian motion $W$ and an independent isotropic $2s$-stable L\'evy process $L$, both $(\cF_t)$-adapted with increments independent of the past.

The coefficients are
\[
b:\R^d\times\cP(\R^d)\to\R^d,\qquad \sigma:\R^d\times\cP(\R^d)\to[0,\infty),\qquad \gamma:\R^d\times\cP(\R^d)\to\R^{d\times d},
\]
and we write $\Sigma:=\sqrt{2\sigma}$. The McKean--Vlasov equation with law-dependent jumps is
\begin{equation}\label{eq:MV-jump}
    \dd X_t = b(X_t,\mu_t)\,\dd t + \Sigma(X_t,\mu_t)\,\dd W_t + \gamma(X_{t-},\mu_t)\,\dd L_t, \qquad \mu_t = \Law(X_t).
\end{equation}
When $\gamma\equiv I_d$, the jump noise is additive and \eqref{eq:MV-jump} becomes
\begin{equation}\label{eq:MV}
  \dd X_t=b(X_t,\mu_t)\,\dd t+\Sigma(X_t,\mu_t)\,\dd W_t+\dd L_t,
  \qquad \mu_t=\Law(X_t).
\end{equation}
A \emph{solution} on $[0,T]$ is a c\`adl\`ag adapted process $X$ satisfying the integral form of the equation with $\mu_t=\Law(X_t)$. When $\sigma\equiv0$ and $\gamma\equiv I_d$, we call the solution a \emph{McKean--Vlasov stable flight}.

We use the following assumptions.

\begin{assumption}\label{ass:H1}
There is $K<\infty$ such that, for all $x,y\in\R^d$ and $\mu,\nu\in\cP(\R^d)$,
\[
  |b(x,\mu)-b(y,\nu)|+|\Sigma(x,\mu)-\Sigma(y,\nu)|\le K\big(|x-y|+\Wr(\mu,\nu)\big).
\]
\end{assumption}

\begin{assumption}\label{ass:H2}
Assumption~\ref{ass:H1} holds, and $b$ and $\sigma$ are bounded.
\end{assumption}

\begin{assumption}\label{ass:Hg}
There is $K<\infty$ such that, for all $x,y\in\R^d$ and $\mu,\nu\in\cP(\R^d)$,
\[
|b(x,\mu)-b(y,\nu)| + |\Sigma(x,\mu)-\Sigma(y,\nu)| + |\gamma(x,\mu)-\gamma(y,\nu)| \le K \bigl( |x-y|+\Wr(\mu,\nu) \bigr).
\]
\end{assumption}

Assumption~\ref{ass:H1} is assumption (H) of the introduction, and Assumption~\ref{ass:Hg} is (H$_\gamma$). Since $\Wr\le1$, Assumption~\ref{ass:H1} implies the linear growth bound $|b(x,\mu)|+|\Sigma(x,\mu)|\le C_0(1+|x|)$,
uniformly in $\mu$.

\begin{example}[Interaction kernels]\label{ex:kernel}
Let $\beta:\R^d\times\R^d\to\R^d$ be bounded and Lipschitz, and set
$\displaystyle b(x,\mu)=\int\beta(x,y)\,\mu(\dd y)$. By Lemma~\ref{lem:Wrho}(c),
Assumption~\ref{ass:H2} holds for $b$. This covers smooth bounded pairwise
interactions $\beta(x,y)=-\nabla V(x-y)$ with $\nabla V$ bounded and
Lipschitz. Unbounded interactions such as $\displaystyle \int(x-y)\,\mu(\dd y)$ are
excluded, which is consistent with Remark~\ref{rem:moments}: the mean need
not exist.
\end{example}

For $\mu\in\cP(\R^d)$ and $f\in C_b^2(\R^d)$, the generators of \eqref{eq:MV} and of \eqref{eq:MV-jump} with frozen law $\mu$ are
\begin{equation}\label{eq:gen-nl}
  \cA_\mu f(x):=b(x,\mu)\cdot\nabla f(x)+\sigma(x,\mu)\Delta f(x)-(-\Delta)^sf(x),
\end{equation}
\begin{equation}\label{eq:gen-jump}
\begin{aligned}
    \cA^\gamma_\mu f(x) &:= b(x,\mu)\cdot\nabla f(x) + \sigma(x,\mu)\Delta f(x) \\
    &\quad + \int_{\R^d\setminus\{0\}} \Bigl[ f\bigl(x+\gamma(x,\mu)y\bigr) - f(x) - \nabla f(x)\cdot\gamma(x,\mu)y\,\one_{\{|y|<1\}} \Bigr] \, \nu_s(\dd y).
\end{aligned}
\end{equation}
Applying them to $e_\xi$, and using \eqref{eq:cos} together with the fact that the odd part of the integrand integrates to zero, one finds their \emph{nonlinear symbols}
\begin{equation}\label{eq:symbol-nl}
  p(x,\xi;\mu):=-e_{-\xi}(x)\,\cA_\mu e_\xi(x)
  =-\iu\, b(x,\mu)\cdot\xi+\sigma(x,\mu)|\xi|^2+|\xi|^{2s},
\end{equation}
\begin{equation}\label{eq:symbol-jump}
    p^\gamma(x,\xi;\mu) := -e_{-\xi}(x)\,\cA^\gamma_\mu e_\xi(x) = -\iu\,b(x,\mu)\cdot\xi + \sigma(x,\mu)|\xi|^2 + |\gamma(x,\mu)^\top\xi|^{2s}.
\end{equation}
In \eqref{eq:symbol-jump}, the law enters the jump part of the symbol.

Finally, let $(X^i_0,W^i,L^i)_{i\ge1}$ be i.i.d.\ copies of $(X_0,W,L)$. The associated interacting particle system is
\begin{equation}\label{eq:particles}
\dd X^{i,N}_t = b(X^{i,N}_t,\mu^N_t)\,\dd t + \Sigma(X^{i,N}_t,\mu^N_t)\,\dd W^i_t + \gamma(X^{i,N}_{t-},\mu^N_{t-})\,\dd L^i_t, \qquad \mu^N_t := \frac{1}{N}\sum_{j=1}^N\delta_{X^{j,N}_t},
\end{equation}
with $X^{i,N}_0=X^i_0$, where $\mu^N_{t-}$ is the empirical measure of the left limits $X^{j,N}_{t-}$, so that the jump integrand is predictable. We denote by $\bar X^i$ the solution of \eqref{eq:MV-jump} driven by $(X^i_0,W^i,L^i)$, and by
\[
\bar\mu^N_t := \frac{1}{N}\sum_{j=1}^N\delta_{\bar X^j_t}
\]
the empirical measure of these independent copies. In Section~\ref{sec:nonlinear}, $\gamma\equiv I_d$.

\section{McKean--Vlasov equations driven by stable noise}
\label{sec:nonlinear}

In this section the jump noise is additive, $\gamma\equiv I_d$, and we study \eqref{eq:MV} under Assumption~\ref{ass:H1}, with the notation of Section~\ref{sec:prelim}.

\begin{theorem}[Strong well-posedness without moments]\label{thm:wellposed}
Let Assumption~\ref{ass:H1} hold, and let $T>0$ and $\mu_0\in\cP(\R^d)$ be
arbitrary. Then \eqref{eq:MV} has a solution on $[0,T]$, and any two
solutions are indistinguishable. Moreover, $t\mapsto\mu_t$ belongs to
$C([0,T];(\cP(\R^d),\Wr))$, and the flow $(\mu_t)$ depends only on $\mu_0$.
\end{theorem}

Let $\mathcal C_T:=C\big([0,T];(\cP(\R^d),\Wr)\big)$, equipped with the
metric
$d_t(m,m'):=\sup_{r\le t}\Wr(m_r,m'_r)$. By Lemma~\ref{lem:Wrho}(a),
$(\mathcal C_T,d_T)$ is complete. For $m\in\mathcal C_T$, let $X^m$ be the
solution of the linear SDE with frozen flow
\begin{equation}\label{eq:frozen}
  X^m_t=X_0+\int_0^tb(X^m_r,m_r)\,\dd r+\int_0^t\Sigma(X^m_r,m_r)\,\dd W_r+L_t .
\end{equation}
The coefficients $(r,x)\mapsto b(x,m_r)$ and $\Sigma(x,m_r)$ are jointly
continuous and $K$-Lipschitz in $x$. Hence the operator
$F(X)_r:=(b(X_r,m_r),\Sigma(X_r,m_r))$ is functional Lipschitz, and
\cite[Thm.~V.7]{Protter2005} yields a unique c\`adl\`ag solution, with no
integrability condition on $X_0$. The solution is the limit of Picard
iterates, each of which is a measurable functional of $(X_0,W,L)$.
Therefore $\Law(X^m)$ depends on $(X_0,W,L)$ only through their (product)
law.

\begin{lemma}[Stability]\label{lem:stability}
Under Assumption~\ref{ass:H1}, there is $C=C(d,K,T)$ such that, for all
$m,m'\in\mathcal C_T$ and $t\le T$,
\[
  \E\Big[\sup_{u\le t}|X^m_u-X^{m'}_u|^2\Big]\le C\int_0^t\Wr(m_r,m'_r)^2\,\dd r .
\]
\end{lemma}

\begin{proof}
Let $D := X^m - X^{m'}$. The term $X_0+L$ is common to both processes, so
\[
D_t = \int_0^t \Delta b_r \, \dd r + \int_0^t \Delta\Sigma_r \, \dd W_r,
\]
where $\Delta b_r := b(X^m_r,m_r) - b(X^{m'}_r,m'_r)$ and $\Delta\Sigma_r$ is defined similarly. Thus $D$ is a continuous semimartingale with $D_0=0$, even though $X^m$ itself may have no moments.

By Assumption~\ref{ass:H1}, $|\Delta b_r|^2 + |\Delta\Sigma_r|^2 \le 2K^2(|D_r|^2+w_r^2)$ with $w_r := \Wr(m_r,m'_r)$. Set $\tau_n := \inf\{t : |D_t| \ge n\}$ and $e_n(t) := \E [ \sup_{u\le t\wedge\tau_n} |D_u|^2 ] \le n^2$. The Cauchy--Schwarz inequality for the $\dd r$-integral and Doob's inequality for the $\dd W$-integral (which brings a factor $d$, as $\Sigma$ is scalar and $W$ is $d$-dimensional) give
\begin{align*}
    e_n(t) &\le 2T \, \E \!\int_0^{t\wedge\tau_n} |\Delta b_r|^2 \dd r + 8d \, \E \!\int_0^{t\wedge\tau_n} |\Delta\Sigma_r|^2 \dd r \\[2ex]
    &\le C_1 \int_0^t \bigl( e_n(r) + w_r^2 \bigr) \dd r,
\end{align*}
with $C_1 := 4K^2(T+4d)$. By Gronwall's lemma, $\displaystyle e_n(t) \le C_1 e^{C_1T} \int_0^t w_r^2 \, \dd r$, and we let $n\to\infty$ by monotone convergence, since $\tau_n \uparrow \infty$ by continuity of $D$.
\end{proof}

\begin{proof}[Proof of Theorem~\ref{thm:wellposed}]

\emph{Step 1: the map $\Phi$.} For $m\in\mathcal C_T$, set $\Phi(m)_t := \Law(X^m_t)$. The process $X^m-L$ is continuous and $L$ is continuous in probability, so $X^m$ is continuous in probability, and by Lemma~\ref{lem:Wrho}(b) and dominated convergence
\[
\Wr\bigl(\Phi(m)_t,\Phi(m)_u\bigr) \le \E \Bigl[ |X^m_t-X^m_u| \wedge 1 \Bigr] \xrightarrow[u\to t]{} 0.
\]
Thus $\Phi$ maps $\mathcal C_T$ into itself and, by the remark following \eqref{eq:frozen}, it depends only on $\mu_0$.

\emph{Step 2: contraction.} Lemma~\ref{lem:Wrho}(b) and Lemma~\ref{lem:stability} give
\begin{align*}
    d_t\bigl(\Phi(m),\Phi(m')\bigr)^2 &\le \E \Biggl[ \sup_{u\le t} |X^m_u-X^{m'}_u|^2 \Biggr] \\[2ex]
    &\le C \int_0^t d_r(m,m')^2 \, \dd r,
\end{align*}
and by induction
\[
d_T\bigl(\Phi^k(m),\Phi^k(m')\bigr)^2 \le \frac{(CT)^k}{k!} d_T(m,m')^2.
\]
For $k$ large, $\Phi^k$ is a strict contraction of the complete space $\mathcal C_T$, so $\Phi$ has a unique fixed point $\mu$, and $X := X^\mu$ solves \eqref{eq:MV}.

\emph{Step 3: uniqueness.} Let $Y$ be a solution and $\nu_t := \Law(Y_t)$. As in Step~1, $\nu\in\mathcal C_T$, and $Y$ solves \eqref{eq:frozen} with $m=\nu$, so $Y=X^\nu$ by uniqueness for \eqref{eq:frozen}. Hence $\nu=\Phi(\nu)$, so $\nu=\mu$ and $Y = X^\mu = X$ up to indistinguishability.
\end{proof}

We now consider the particle system \eqref{eq:particles} with $\gamma\equiv I_d$.

\begin{theorem}\label{thm:chaos}
Let Assumption~\ref{ass:H1} hold and let $\mu_0\in\cP(\R^d)$ be arbitrary. Then \eqref{eq:particles} is well posed, and
\[
\lim_{N\to\infty} \E \Biggl[ \sup_{t\le T} |X^{1,N}_t-\bar X^1_t|^2 \Biggr] = 0.
\]

In particular, for every fixed $k$, the law of $(X^{1,N},\dots,X^{k,N})$ converges to the product measure $\Law(\bar X^1)^{\otimes k}$, the law of $k$ independent copies of $\bar X^1$, on the Skorokhod space $D([0,T];\R^d)^k$.
\end{theorem}

\begin{proof}
\emph{Well-posedness.} For $\mathbf x,\mathbf y\in(\R^d)^N$, Lemma~\ref{lem:Wrho}(b) applied to the uniform coupling of the indices gives
\[
\Wr \Biggl( \frac{1}{N}\sum_j\delta_{x^j}, \, \frac{1}{N}\sum_j\delta_{y^j} \Biggr) \le \frac{1}{N}\sum_j |x^j-y^j|.
\]
So \eqref{eq:particles} is an SDE in $\R^{Nd}$ with Lipschitz coefficients and additive jump noise, which is well posed by \cite[Thm.~V.7]{Protter2005}. Its solution is a measurable functional of $(X^i_0,W^i,L^i)_{i\le N}$ that commutes with permutations of the indices, so the system is exchangeable.

\emph{Coupling estimate.} Let $D^i := X^{i,N}-\bar X^i$, and recall the empirical measure $\bar\mu^N_t$ of Section~\ref{sec:prelim}. As before, $X_0^i$ and $L^i$ cancel, so $D^i$ is a continuous semimartingale. By the triangle inequality and the coupling bound,
\[
|b(X^{i,N}_r,\mu^N_r)-b(\bar X^i_r,\mu_r)| \le K \Biggl( |D^i_r| + \frac{1}{N}\sum_j|D^j_r| + \Wr(\bar\mu^N_r,\mu_r) \Biggr),
\]
and similarly for $\Sigma$. Let
\[
\tau_n := \inf \Bigl\{ t : \sum_j|D^j_t| \ge n \Bigr\},
\]
which is symmetric in the indices, and $e_n(t) := \E [ \sup_{u\le t\wedge\tau_n} |D^1_u|^2 ]$. By exchangeability,
\[
\E \Biggl[ \sup_{u\le t\wedge\tau_n} \frac{1}{N}\sum_j|D^j_u|^2 \Biggr] \le \frac{1}{N}\sum_j \E \Biggl[ \sup_{u\le t\wedge\tau_n}|D^j_u|^2 \Biggr] = e_n(t),
\]
so the argument of Lemma~\ref{lem:stability} gives
\[
e_n(t) \le C \int_0^t \Bigl( e_n(r) + \E \bigl[ \Wr(\bar\mu^N_r,\mu_r)^2 \bigr] \Bigr) \, \dd r,
\]
with $C=C(d,K,T)$. Gronwall's lemma and $n\to\infty$ give
\[
\E \Biggl[ \sup_{t\le T}|D^1_t|^2 \Biggr] \le C e^{CT} \int_0^T \E \bigl[ \Wr(\bar\mu^N_r,\mu_r)^2 \bigr] \, \dd r.
\]

\emph{Law of large numbers.} For each $r$, the variables $\bar X^j_r$ are i.i.d.\ with law $\mu_r$. By Varadarajan's theorem \cite{Varadarajan1958}, $\bar\mu^N_r\to\mu_r$ weakly almost surely, hence $\Wr(\bar\mu^N_r,\mu_r)\to0$ almost surely by Lemma~\ref{lem:Wrho}(a). Since $\Wr\le1$, the right-hand side tends to $0$ by dominated convergence, first in $\omega$ and then in $r$.

\emph{Chaos.} The limits $\bar X^i$ are i.i.d., and convergence in $L^2$ for the uniform distance implies convergence in law for the Skorokhod topology, which gives the statement for $k$ particles.
\end{proof}
\begin{remark}[Rates]
The proof gives the rate of convergence for the integral $\displaystyle \int_0^T \E [ \Wr(\bar\mu^N_r,\mu_r)^2 ] \, \dd r$.
For general coefficients and without moments, this convergence can be arbitrarily slow. Since $\Wr\le W_1$, rates can be obtained from \cite{FournierGuillin2015} when $s>\tfrac12$ and $\mu_r$ has a finite moment of some order $q\in(1,2s)$. For interactions of kernel type, the following theorem gives the rate $N^{-1/2}$ with no moment assumption at all.
\end{remark}

\begin{theorem}[Chaos at rate $N^{-1/2}$ for kernel interactions]\label{thm:chaos-rate}
Let $\beta:\R^d\times\R^d\to\R^d$ and $\varsigma:\R^d\times\R^d\to[0,\infty)$ be bounded by $M$ and $\ell$-Lipschitz on $\R^d\times\R^d$. We define the coefficients
\begin{align*}
    b(x,\mu) &= \int \beta(x,y) \, \mu(\dd y), \\[1.5ex]
    \Sigma(x,\mu) &= \int \varsigma(x,y) \, \mu(\dd y),
\end{align*}
that is, $\sigma=\Sigma^2/2$. Then Assumption~\ref{ass:H2} holds, and, with the notation of Theorem~\ref{thm:chaos}, there is a constant $C=C(d,M,\ell,T)$, independent of $\mu_0$ and $s$, such that for all $N\ge1$ and $i\le N$
\[
\E \Biggl[ \sup_{t\le T} |X^{i,N}_t-\bar X^i_t|^2 \Biggr] \le \frac{C}{N}.
\]
\end{theorem}

\begin{proof}
Assumption~\ref{ass:H2} follows from Lemma~\ref{lem:Wrho}(c), applied componentwise and uniformly in $x$, and from the Lipschitz property in $x$. We keep the notation $D^i=X^{i,N}-\bar X^i$ and $\bar\mu^N_r$ of the proof of Theorem~\ref{thm:chaos}, and write
\[
b(X^{i,N}_r,\mu^N_r) - b(\bar X^i_r,\mu_r) = \Bigl[ b(X^{i,N}_r,\mu^N_r) - b(\bar X^i_r,\bar\mu^N_r) \Bigr] + R^i_r,
\]
with
\[
R^i_r := \frac{1}{N}\sum_{j=1}^N \xi^{ij}_r, \qquad \xi^{ij}_r := \beta(\bar X^i_r,\bar X^j_r) - \int \beta(\bar X^i_r,y) \, \mu_r(\dd y).
\]

\emph{First term.} Since $\beta$ is $\ell$-Lipschitz,
\begin{align*}
    \bigl| b(X^{i,N}_r,\mu^N_r) - b(\bar X^i_r,\bar\mu^N_r) \bigr| &\le \frac{1}{N}\sum_j \bigl| \beta(X^{i,N}_r,X^{j,N}_r) - \beta(\bar X^i_r,\bar X^j_r) \bigr| \\[2ex]
    &\le \ell \Biggl( |D^i_r| + \frac{1}{N}\sum_j|D^j_r| \Biggr).
\end{align*}

\emph{Fluctuation term.} The variables $\bar X^1_r,\dots,\bar X^N_r$ are i.i.d.\ with law $\mu_r$, since they are driven by independent noises. Let $j\ne k$; one of them, say $k$, differs from $i$. Conditionally on $(\bar X^i_r,\bar X^j_r)$, the variable $\bar X^k_r$ still has law $\mu_r$, so $\E [ \xi^{ik}_r \,|\, \bar X^i_r,\bar X^j_r ] = 0$, and since $\xi^{ij}_r$ is a function of $(\bar X^i_r,\bar X^j_r)$, we get $\E[\xi^{ij}_r\cdot\xi^{ik}_r]=0$ (this includes the case $j=i$). Only the diagonal terms remain:
\[
\E \bigl[ |R^i_r|^2 \bigr] = \frac{1}{N^2}\sum_{j=1}^N \E \bigl[ |\xi^{ij}_r|^2 \bigr] \le \frac{4M^2}{N}.
\]
The same holds for $\Sigma$, with $\varsigma$ in place of $\beta$ and a remainder $\tilde R^i_r$ such that $\E[|\tilde R^i_r|^2] \le 4M^2/N$.

\emph{Gronwall.} With $\tau_n$ and $e_n$ as in the proof of Theorem~\ref{thm:chaos}, the Cauchy--Schwarz and Doob inequalities and exchangeability give
\begin{align*}
    e_n(t) &\le C_1 \int_0^t \Bigl( e_n(r) + \E \bigl[ |R^1_r|^2 \bigr] + \E \bigl[ |\tilde R^1_r|^2 \bigr] \Bigr) \, \dd r \\[2ex]
    &\le C_1 \int_0^t e_n(r) \, \dd r + \frac{8C_1 M^2 T}{N},
\end{align*}
where $C_1$ depends only on $d,\ell,T$. Gronwall's lemma and $n\to\infty$ give the claim with $C = 8C_1 M^2 T e^{C_1 T}$. The initial law and the stability index do not enter the constants, because the initial conditions and the stable noises cancel in $D^i$.
\end{proof}

\begin{remark}
Theorem~\ref{thm:chaos-rate} is the analogue of McKean's theorem \cite[Thm.~3.1]{ChaintronDiez2022} for additive stable noise; the rate $N^{-1/2}$ in $L^2$ is that of the law of large numbers for the fluctuation term $R^i$.
\end{remark}

We end this section with further properties of the flow: McKean's nonlinear semigroup, the martingale problem with its carr\'e du champ, the nonlinear symbol, and explicit laws when the coefficients depend on the measure only.

\begin{proposition}[McKean's nonlinear semigroup]\label{prop:nl-semigroup}
Under Assumption~\ref{ass:H1}, for $\mu_0\in\cP(\R^d)$ let $T_t(\mu_0) := \mu_t = \Law(X_t)$, where $X$ is the solution of \eqref{eq:MV}.
Then $t\mapsto T_t(\mu_0)$ is $\Wr$-continuous, and the family $(T_t)_{t\ge 0}$ forms a semigroup satisfying
\begin{align*}
    T_0 &= \mathrm{Id}, \\[1.5ex]
    T_{t+u} &= T_t \circ T_u \qquad \text{for all } t,u \ge 0.
\end{align*}
Moreover, given the flow $(\mu_t)$, the process $X$ is a time-inhomogeneous Markov process with time-dependent generator $\cA_{\mu_t}$.
\end{proposition}

\begin{proof}
Fix $u\ge0$ and set $Y_t := X_{u+t}$, $W'_t := W_{u+t} - W_u$, $L'_t := L_{u+t} - L_u$ and $\mathcal G_t := \cF_{u+t}$. Then $W'$ is a $(\mathcal G_t)$-Brownian motion, $L'$ is an independent $(\mathcal G_t)$-stable process, and $Y_0=X_u$ is $\mathcal G_0$-measurable. So $Y$ solves \eqref{eq:MV} driven by $(W',L')$, and its flow is $\Law(Y_t) = \mu_{u+t}$. By Theorem~\ref{thm:wellposed}, this flow depends only on $\mu_u$, which is the semigroup property $\mu_{u+t} = T_t(\mu_u)$. Continuity in $t$ was proved in Theorem~\ref{thm:wellposed}. Finally, $X$ is a Markov process because it solves the linear equation \eqref{eq:frozen} with the deterministic coefficients $(r,x) \mapsto ( b(x,\mu_r), \Sigma(x,\mu_r) )$.
\end{proof}

\begin{remark}
$(T_t)$ is a semigroup of nonlinear maps on $\cP(\R^d)$, which is not
a linear space, so it does not come from Hille--Yosida theory; linear
semigroup theory only describes the frozen dynamics. This is what
``nonlinear in the sense of McKean'' \cite{McKean1966} means.
\end{remark}

\begin{theorem}\label{thm:MP}
Let Assumption~\ref{ass:H1} hold and let $X$ be the solution of \eqref{eq:MV}. For every $\varphi\in C_c^\infty(\R^d)$, the process
\begin{equation}\label{eq:MP}
    M^\varphi_t := \varphi(X_t) - \varphi(X_0) - \int_0^t \cA_{\mu_r}\varphi(X_r) \, \dd r
\end{equation}
is a square-integrable martingale with predictable quadratic variation given by
\[
\langle M^\varphi\rangle_t = \int_0^t \Gamma_{\mu_r}(\varphi,\varphi)(X_r) \, \dd r.
\]

Here the \emph{nonlinear carr\'e du champ} is
\begin{equation}\label{eq:carre-nl}
\begin{aligned}
    \Gamma_\mu(f,g)(x) &:= \cA_\mu(fg) - f\cA_\mu g - g\cA_\mu f \\[1ex]
    &= 2\sigma(x,\mu)\nabla f\cdot\nabla g + \int \bigl(f(x+y)-f(x)\bigr)\bigl(g(x+y)-g(x)\bigr) \, \nu_s(\dd y).
\end{aligned}
\end{equation}

In particular, $(\mu_t)$ is a weak solution of the nonlinear Fokker--Planck equation $\partial_t\mu_t = \cA_{\mu_t}^*\mu_t$: the map $t \mapsto \ip{\mu_t}{\varphi}$ belongs to $C^1$ and satisfies
\begin{equation}\label{eq:FP-nl}
    \ip{\mu_t}{\varphi} = \ip{\mu_0}{\varphi} + \int_0^t \ip{\mu_r}{\cA_{\mu_r}\varphi} \, \dd r.
\end{equation}

Under Assumption~\ref{ass:H2}, all of the above holds for every $\varphi\in C_b^2(\R^d)$, including complex-valued $\varphi$.
\end{theorem}

\begin{proof}
Given the flow $(\mu_r)$, the process $X$ solves the linear equation \eqref{eq:frozen}, and the It\^o computation of Proposition~\ref{prop:ito} applies with $b(\cdot)$ and $\sigma(\cdot)$ replaced by $b(\cdot,\mu_r)$ and $\sigma(\cdot,\mu_r)$. Formula \eqref{eq:carre-nl} is obtained as \eqref{eq:carre}: for fixed $\mu$, the operator $b(\cdot,\mu)\cdot\nabla$ is a derivation, so the dependence of $b$ on the measure does not enter $\Gamma_\mu$.

It remains to check integrability. Let $\varphi\in C_c^\infty$ be supported in a ball $B$. By the linear growth bound and $\Wr \le 1$, $\sup_{r\le T,\,x\in B} ( |b(x,\mu_r)| + \sigma(x,\mu_r) ) < \infty$, so $\cA_{\mu_r}\varphi$ is bounded on $[0,T]\times\R^d$. It is also jointly continuous in $(r,x)$, since $r\mapsto\mu_r$ is $\Wr$-continuous and $b$, $\sigma$ are continuous. The Brownian integrand is bounded, and the jump integrand is handled as in Proposition~\ref{prop:ito}.

Taking expectations in \eqref{eq:MP} gives \eqref{eq:FP-nl}. The map $r \mapsto \E [ \cA_{\mu_r}\varphi(X_r) ]$ is continuous, by continuity in probability of $X$ and dominated convergence, so $t\mapsto\ip{\mu_t}{\varphi}$ is $C^1$. Under Assumption~\ref{ass:H2}, the same bounds hold for every $\varphi\in C_b^2$.
\end{proof}

\begin{remark}[Where the measure derivative enters]
The dependence of the coefficients on the measure does not appear in $\Gamma_\mu$; it appears only when one differentiates functionals of the flow. Formally, for a smooth $U:\cP(\R^d)\to\R$ with linear functional derivative $\displaystyle\frac{\delta U}{\delta\mu}(\mu,y)$, equation \eqref{eq:FP-nl} and the chain rule \cite[Ch.~5]{CarmonaDelarue2018} give
\[
\frac{\dd}{\dd t}U(\mu_t) = \int \cA_{\mu_t} \Biggl[ \frac{\delta U}{\delta\mu}(\mu_t,\cdot) \Biggr](x) \, \mu_t(\dd x),
\]
and the Lions derivative is $\displaystyle \partial_\mu U(\mu)(y) = \nabla_y \frac{\delta U}{\delta\mu}(\mu,y)$. These objects govern the linearisation of $T_t$ and the master equation, which we do not pursue here.
\end{remark}

\begin{theorem}\label{thm:nl-symbol}
Let Assumption~\ref{ass:H2} hold and let $\chi_t := \chi_{\mu_t}$. Then, for all $\xi\in\R^d$, we have
\begin{equation}\label{eq:char-nl}
    \frac{\dd}{\dd t}\chi_t(\xi) = -\int_{\R^d} p(x,\xi;\mu_t) \, e^{\iu\xi\cdot x} \, \mu_t(\dd x).
\end{equation}

Moreover, if $X^x$ denotes the solution of \eqref{eq:MV} with $\mu_0 = \delta_x$, then
\begin{equation}\label{eq:schilling-nl}
    p(x,\xi;\delta_x) = -\lim_{t\downarrow0} \frac{1}{t} \Biggl( \E \Bigl[ e^{\iu\xi\cdot(X^x_t-x)} \Bigr] - 1 \Biggr).
\end{equation}
\end{theorem}

\begin{proof}
Formula \eqref{eq:char-nl} is \eqref{eq:FP-nl} with $\varphi = e_\xi$, together with \eqref{eq:symbol-nl}. For \eqref{eq:schilling-nl}, \eqref{eq:FP-nl} with $\mu_0 = \delta_x$ gives
\[
\frac{1}{t} \Biggl( \E \Bigl[ e^{\iu\xi\cdot(X^x_t-x)} \Bigr] - 1 \Biggr) = -\frac{1}{t} \int_0^t \E \Bigl[ p(X^x_r,\xi;\mu_r) e^{\iu\xi\cdot(X^x_r-x)} \Bigr] \, \dd r.
\]
As $r \downarrow 0$, $X^x_r \to x$ almost surely and $\mu_r \to \delta_x$ in $\Wr$. Since $p(\cdot,\xi;\cdot)$ is bounded and continuous under Assumption~\ref{ass:H2}, the integrand is right-continuous at $r=0$ by dominated convergence, and the limit follows.
\end{proof}

Formula \eqref{eq:schilling-nl} is the nonlinear analogue of \eqref{eq:schilling}: at small times, a particle started at $x$ sees the symbol of the frozen generator $\cA_{\delta_x}$.

\begin{corollary}[Measure-dependent coefficients: explicit laws]\label{cor:explicit}
Let Assumption~\ref{ass:H1} hold with $b(x,\mu)=b(\mu)$ and $\sigma(x,\mu)=\sigma(\mu)$ independent of $x$, and let
\[
B_t := \int_0^t b(\mu_r) \, \dd r \qquad \text{and} \qquad S_t := \int_0^t \sigma(\mu_r) \, \dd r.
\]
Then
\begin{equation}\label{eq:explicit}
    \chi_t(\xi) = \chi_0(\xi) \exp\Biggl(-\int_0^t p(\xi;\mu_r) \, \dd r \Biggr) = \chi_0(\xi) \, e^{\iu\xi\cdot B_t - S_t|\xi|^2 - t|\xi|^{2s}},
\end{equation}
that is, $\mu_t = \mu_0 * \delta_{B_t} * \mathcal N(0,2S_tI_d) * \pi_t$. Moreover, $(B,S)$ is the unique solution of the $(d+1)$-dimensional ODE system
\begin{equation}\label{eq:ODE}
\begin{cases}
    \dot B_t = b\bigl(G_t(B_t,S_t)\bigr), \\[1.5ex]
    \dot S_t = \sigma\bigl(G_t(B_t,S_t)\bigr), \\[1.5ex]
    (B_0,S_0) = (0,0),
\end{cases}
\end{equation}
where $G_t(B,S) := \mu_0 * \delta_B * \mathcal N(0,2SI_d) * \pi_t$.
\end{corollary}

\begin{proof}
Here $\displaystyle X_t = X_0 + B_t + \int_0^t \Sigma(\mu_r) \, \dd W_r + L_t$. The Wiener integral has a deterministic integrand, so it is Gaussian with law
\[
\mathcal N \Biggl( 0, \int_0^t \Sigma(\mu_r)^2 \dd r \, I_d \Biggr) = \mathcal N(0, 2S_t I_d).
\]
The three random terms are independent, which gives \eqref{eq:explicit}, and in particular $(B,S)$ solves \eqref{eq:ODE}.

Conversely, let $(\tilde B,\tilde S)$ be a $C^1$ solution of \eqref{eq:ODE}; then $\dot{\tilde S} \ge 0$, and we set
\[
\tilde X_t := X_0 + \tilde B_t + \int_0^t \sqrt{2\dot{\tilde S}_r} \, \dd W_r + L_t.
\]
The law of $\tilde X_t$ is $G_t(\tilde B_t,\tilde S_t) =: \tilde\mu_t$, and $\dot{\tilde B}_t = b(\tilde\mu_t)$, $\sqrt{2\dot{\tilde S}_t} = \Sigma(\tilde\mu_t)$, so $\tilde X$ solves \eqref{eq:MV}. By the uniqueness part of Theorem~\ref{thm:wellposed}, $\tilde\mu = \mu$, hence $(\tilde B,\tilde S) = (B,S)$.
\end{proof}

\begin{example}[A mean-field stable flight]\label{ex:flight}
Let $d=1$, $\kappa>0$, $\sigma\equiv0$, $\mu_0=\delta_{x_0}$ and $\displaystyle b(\mu) = -\kappa \int \tanh(y) \, \mu(\dd y)$. By Example~\ref{ex:kernel}, $b$ satisfies Assumption~\ref{ass:H2}. By Corollary~\ref{cor:explicit}, $\mu_t = \delta_{x_0+B_t} * \pi_t$, where
\[
\dot B_t = -\kappa \, g_t(x_0+B_t), \qquad g_t(z) := \E \bigl[ \tanh(z+L_t) \bigr].
\]
Since $\tanh$ is odd and increasing and $L_t$ is symmetric, $g_t$ is odd and non-decreasing, so the centre $x_0+B_t$ moves monotonically towards $0$. This makes sense for every $s\in(0,1)$, including $s\le\tfrac12$, where $\E|X_t| = \infty$ and there is no mean position.
\end{example}

\section{Law-dependent jumps: a critical threshold at \texorpdfstring{$s=\frac12$}{s=1/2}}
\label{sec:jumps}

We now let the size of the jumps depend on the position and on the law, and study \eqref{eq:MV-jump} under Assumption~\ref{ass:Hg}. Its generator is \eqref{eq:gen-jump} and its symbol \eqref{eq:symbol-jump}.

When $s>\frac12$, the stable process $L$ has finite moments of every order $p<2s$, in particular of some order $p>1$. The following estimate replaces the $L^2$ estimates of Section~\ref{sec:nonlinear}.

\begin{lemma}\label{lem:Lp}
Let $s\in(\frac12,1)$ and $p\in[1,2s)$. There is a constant $C_*$, depending only on $d$, $s$ and $p$, with the following property. Let $t_0<t_1$, let $\beta$ be a progressively measurable $\R^d$-valued process, $\theta$ a progressively measurable real process and $g$ a predictable $\R^{d\times d}$-valued process, all locally bounded (for $\beta$ and $\theta$, it suffices that this holds $\dd r\otimes\dd\Prob$-almost everywhere). We consider the process
\[
D_t = D_{t_0} + \int_{t_0}^t \beta_r\,\dd r + \int_{t_0}^t \theta_r\,\dd W_r + \int_{t_0}^t g_r\,\dd L_r, \qquad t\in[t_0,t_1],
\]
with $D_{t_0}$ being $\cF_{t_0}$-measurable. Then we have the maximal inequality
\begin{align*}
    \E \Biggl[ \sup_{t_0\le t\le t_1}|D_t|^p \Biggr] &\le C_* \Biggl( \E \bigl[ |D_{t_0}|^p \bigr] + \E \Biggl[ \biggl( \int_{t_0}^{t_1}|\beta_r|\dd r \biggr)^p \Biggr] \\
    &\quad + \E \Biggl[ \biggl( \int_{t_0}^{t_1}\bigl(\theta_r^2+|g_r|^2\bigr)\dd r \biggr)^{p/2} \Biggr] + \E \Biggl[ \int_{t_0}^{t_1}|g_r|^p\,\dd r \Biggr] \Biggr).
\end{align*}
\end{lemma}

\begin{proof}
By \eqref{eq:LI}, and since for $s>\frac12$
\[
\int_{\{|y|\ge1\}} |y| \, \nu_s(\dd y) < \infty \qquad \text{and} \qquad \int_{\{|y|\ge1\}} y \, \nu_s(\dd y) = 0,
\]
the jump integral splits as $\displaystyle \int_{t_0}^t g_r\,\dd L_r = M^1_t + M^2_t$, where
\[
M^1_t := \int_{t_0}^t \!\!\int_{\{|y|<1\}} g_r y\,\tilde N(\dd r,\dd y), \qquad M^2_t := \int_{t_0}^t \!\!\int_{\{|y|\ge1\}} g_r y\,\tilde N(\dd r,\dd y)
\]
are local martingales.

It suffices to treat bounded $\beta$, $\theta$ and $g$. Indeed, they are locally bounded, so there are stopping times $\sigma_k\uparrow\infty$ such that the integrands multiplied by $\one_{\{r\le\sigma_k\}}$ are bounded; the right-hand side of the inequality for $D^{\sigma_k}$ is at most that for $D$, and its left-hand side increases to $\E[\sup_{t_0\le t\le t_1}|D_t|^p]$ by monotone convergence. So let $\beta$, $\theta$ and $g$ be bounded; then $M^1$ is a square-integrable martingale. We bound the five terms of $D$ separately, using $(a_1+\dots+a_5)^p \le 5^{p-1}(a_1^p+\dots+a_5^p)$.

The drift term is at most $\displaystyle \Bigl( \int_{t_0}^{t_1} |\beta_r| \, \dd r \Bigr)^p$, and by the Burkholder--Davis--Gundy inequality \cite[Thm.~IV.48]{Protter2005} the Brownian term is at most
\[
C \, \E \Biggl[ \Biggl( \int_{t_0}^{t_1} \theta_r^2\,\dd r \Biggr)^{p/2} \Biggr],
\]
with $C$ depending on $p$ and $d$.

For $M^1$, the Burkholder--Davis--Gundy inequality gives $\E[\sup|M^1|^p] \le C\,\E\bigl[ [M^1]_{t_1}^{p/2} \bigr]$, and the predictable compensator of $[M^1]$ satisfies
\[
\langle M^1\rangle_t = \int_{t_0}^t \int_{\{|y|<1\}} |g_r y|^2 \nu_s(\dd y)\dd r \le c_1 \int_{t_0}^t |g_r|^2 \dd r, \qquad c_1 := \int_{\{|y|<1\}} |y|^2 \, \nu_s(\dd y).
\]
We use Lenglart's domination inequality \cite[Corollaire~II]{Lenglart1977} in the following form: if $A$ is an adapted increasing c\`adl\`ag process and $B$ a predictable increasing process with $\E A_\tau\le\E B_\tau$ for every bounded stopping time $\tau$, then $\E[A_{t_1}^k]\le\frac{2-k}{1-k}\,\E[B_{t_1}^k]$ for every $k\in(0,1)$. Since $M^1$ is square-integrable, $[M^1]-\langle M^1\rangle$ is a martingale, so $\E[M^1]_\tau=\E\langle M^1\rangle_\tau$ for every bounded stopping time $\tau$. Applying the inequality with $A=[M^1]$, $B=\langle M^1\rangle$ and $k=p/2<1$ gives
\[
\E \bigl[ [M^1]_{t_1}^{p/2} \bigr] \le \frac{4-p}{2-p} \, \E \bigl[ \langle M^1\rangle_{t_1}^{p/2} \bigr].
\]

For $M^2$, note first that
\[
\E\int_{t_0}^{t_1}\!\!\int_{\{|y|\ge1\}}|g_ry|\,\nu_s(\dd y)\,\dd r<\infty,
\]
since $g$ is bounded and $\int_{\{|y|\ge1\}}|y|\,\nu_s(\dd y)<\infty$; hence $M^2$ is a martingale. We apply the Burkholder--Davis--Gundy inequality, valid for every $p\ge1$ and every local martingale, and then $(\sum_ia_i)^{p/2} \le \sum_ia_i^{p/2}$, which holds since $p/2\le1$:
\begin{align*}
    \E \Biggl[ \sup |M^2|^p \Biggr] &\le C \, \E \Biggl[ \Biggl( \sum_{t_0<r\le t_1} |g_r\Delta L_r|^2 \one_{\{|\Delta L_r|\ge1\}} \Biggr)^{p/2} \Biggr] \\[2ex]
    &\le C \, \E \Biggl[ \sum_{t_0<r\le t_1} |g_r|^p |\Delta L_r|^p \one_{\{|\Delta L_r|\ge1\}} \Biggr] \\[2ex]
    &= C \, c_p \, \E \Biggl[ \int_{t_0}^{t_1} |g_r|^p\,\dd r \Biggr],
\end{align*}
where $\displaystyle c_p := \int_{\{|y|\ge1\}} |y|^p \, \nu_s(\dd y)$ is finite because $p<2s$. The last equality is the compensation formula for the nonnegative predictable integrand $|g_r|^p|y|^p\one_{\{|y|\ge1\}}$.
\end{proof}

\begin{theorem}\label{thm:jump-wp}
Let $s\in(\frac12,1)$ and let Assumption~\ref{ass:Hg} hold. For every $\mu_0\in\cP(\R^d)$ and every $T>0$, the following hold.

\emph{Strong well-posedness.} Equation \eqref{eq:MV-jump} with $\Law(X_0)=\mu_0$ has a strong solution on $[0,T]$, and any two solutions are indistinguishable. The flow $t\mapsto\mu_t$ is $\Wr$-continuous, depends only on $\mu_0$, and $T_t(\mu_0):=\mu_t$ defines a semigroup of nonlinear maps on $\cP(\R^d)$.

\emph{Fokker--Planck formulation.} If moreover $b$, $\Sigma$ and $\gamma$ are bounded, then for every $\varphi\in C_b^2(\R^d)$, the flow satisfies
\[
\ip{\mu_t}{\varphi} = \ip{\mu_0}{\varphi} + \int_0^t \ip{\mu_r}{\cA^\gamma_{\mu_r}\varphi}\,\dd r.
\]
\end{theorem}

\begin{proof}
Fix $p\in(1,2s)$, let $C_*$ be the constant of Lemma~\ref{lem:Lp}, and choose $T_0\in(0,1]$, depending only on $d,s,p,K$, such that
\begin{equation}\label{eq:T0}
    \vartheta := C_* \, 4^{p}K^p \bigl( T_0^p+2T_0^{p/2}+T_0 \bigr) \le \tfrac14.
\end{equation}

\emph{Step 1: the frozen equation on a short interval.} Let $t_0\ge0$, $t_1:=t_0+T_0$, let $\xi$ be $\cF_{t_0}$-measurable with law $\nu$, and let $m\in C([t_0,t_1];(\cP(\R^d),\Wr))$. We consider the frozen equation
\[
X^m_t = \xi + \int_{t_0}^t b(X^m_r,m_r)\,\dd r + \int_{t_0}^t \Sigma(X^m_r,m_r)\,\dd W_r + \int_{t_0}^t \gamma(X^m_{r-},m_r)\,\dd L_r.
\]
This equation has functional Lipschitz coefficients and is driven by the semimartingale $(t,W,L)$, so it has a unique c\`adl\`ag solution \cite[Thm.~V.7]{Protter2005}, with no integrability condition on $\xi$. As in Section~\ref{sec:nonlinear}, its law depends only on $\nu$ and $m$. Since $L$ has no fixed jump times, $X^m$ is continuous in probability, so $\Phi(m)_t := \Law(X^m_t)$ defines a map $\Phi$ of $C([t_0,t_1];(\cP(\R^d),\Wr))$ into itself.

\emph{Step 2: contraction.} Let $m,m'$ be two flows, $D := X^m - X^{m'}$, $w_r := \Wr(m_r,m'_r)$ and $\bar w := \sup_{[t_0,t_1]} w$. Then $D_{t_0}=0$, and $D$ has the form of Lemma~\ref{lem:Lp} with $\beta_r = b(X^m_r,m_r)-b(X^{m'}_r,m'_r)$ and similar expressions for $\theta$ and $g$, each bounded for almost every $r$ by $K ( |D_{r-}| + w_r )$. Let $\tau_n := \inf\{t\ge t_0 : |D_t|\ge n\} \wedge t_1$ and $a_n := \E [ \sup_{t_0\le t\le\tau_n} |D_t|^p ]$. The stopped process $D^{\tau_n}$ has the same form, with integrands multiplied by $\one_{\{r\le\tau_n\}}$ and bounded by $K(n+1)$, so $a_n<\infty$.

For $r\le\tau_n$ we have $|D_{r-}| \le \sup_{t_0\le t\le\tau_n}|D_t|$, hence, by $(a+b)^p\le 2^{p-1}(a^p+b^p)$,
\[
|\beta_r|^p,\ |\theta_r|^p,\ |g_r|^p \le 2^{p-1}K^p \Bigl( \sup_{t_0\le t\le\tau_n}|D_t|^p + \bar w^p \Bigr).
\]
Using also $(a+b)^{p/2} \le a^{p/2}+b^{p/2}$, Lemma~\ref{lem:Lp} gives
\[
a_n \le C_* \, 4^{p}K^p \bigl( T_0^p+2T_0^{p/2}+T_0 \bigr) \bigl( a_n+\bar w^p \bigr) = \vartheta \, (a_n+\bar w^p).
\]
By \eqref{eq:T0}, $\vartheta \le \frac{1}{4}$, so $a_n \le \frac{1}{3}\bar w^p$, and letting $n\to\infty$ we get $\E [ \sup_{[t_0,t_1]}|D|^p ] \le \frac{1}{3}\bar w^p$ by monotone convergence. With Lemma~\ref{lem:Wrho}(b) and Jensen's inequality ($p\ge1$),
\[
\Wr\bigl(\Phi(m)_t,\Phi(m')_t\bigr) \le \E \Bigl[ |D_t|\wedge1 \Bigr] \le \Bigl( \E \bigl[ |D_t|^p \bigr] \Bigr)^{1/p} \le 3^{-1/p} \, \bar w,
\]
so $\Phi$ is a strict contraction for the uniform distance. Its unique fixed point gives a solution of \eqref{eq:MV-jump} on $[t_0,t_1]$ starting from $\xi$.

\emph{Step 3: concatenation and uniqueness.} Since $T_0$ does not depend on the initial law, we can apply Step~2 successively on $[kT_0,(k+1)T_0]$, starting at time $kT_0$ from the value of the solution already constructed. This gives a solution on $[0,T]$. If $Y$ is another solution, its flow $\nu_t:=\Law(Y_t)$ is $\Wr$-continuous, since $\Delta Y_t=\gamma(Y_{t-},\nu_t)\Delta L_t=0$ almost surely at each fixed $t$; on $[0,T_0]$ it is therefore a fixed point of $\Phi$, so it coincides with that of $X$, and $Y=X$ on $[0,T_0]$ by uniqueness for the frozen equation. By induction, $Y=X$ on $[0,T]$. The remaining statements for the first part follow as in Theorem~\ref{thm:wellposed} and Proposition~\ref{prop:nl-semigroup}.

For the Fokker--Planck equation, we apply It\^o's formula as in Proposition~\ref{prop:ito}, with the jump $y$ replaced by $\gamma(X_{r-},\mu_r)y$. Since $b$, $\Sigma$ and $\gamma$ are bounded, $\cA^\gamma_{\mu_r}\varphi(X_r)$ is bounded, the Brownian integral is a martingale, and the compensated Poisson integral is a square-integrable martingale, because
\[
\int |\varphi(x+\gamma y)-\varphi(x)|^2 \nu_s(\dd y) \le \int \min \bigl( \norm{\nabla\varphi}_\infty^2|\gamma|^2|y|^2, \, 4\norm\varphi_\infty^2 \bigr) \, \nu_s(\dd y) < \infty.
\]
Taking expectations gives the equation.
\end{proof}

\begin{remark}\label{rem:where-s}
The assumption $s>\frac12$ is used twice: to have
\[
\int_{\{|y|\ge1\}} |y|^p \, \nu_s(\dd y) < \infty
\]
for some $p\ge1$, and to get the \emph{linear} bound $\Wr(\Phi(m)_t,\Phi(m')_t) \le (\E|D_t|^p)^{1/p}$.
With moments of order $p<1$ only, one gets $\Wr \le \E|D_t|^p \lesssim \bar w^{\,p}$, which is not a contraction, and Proposition~\ref{prop:counter} shows that this is not a weakness of the method.
\end{remark}

We now consider the particle system \eqref{eq:particles} with a general $\gamma$.

\begin{theorem}\label{thm:jump-chaos}
Let $s\in(\frac12,1)$, $p\in[1,2s)$, let Assumption~\ref{ass:Hg} hold and let $\mu_0\in\cP(\R^d)$. Then, for every $T>0$,
\[
\lim_{N\to\infty} \E \Biggl[ \sup_{t\le T} |X^{i,N}_t-\bar X^i_t|^p \Biggr] = 0.
\]
Assume moreover that
\[
b(x,\mu) = \int \beta(x,y)\,\mu(\dd y), \qquad
\Sigma(x,\mu) = \int \varsigma(x,y)\,\mu(\dd y), \qquad
\gamma(x,\mu) = \int \kappa(x,y)\,\mu(\dd y),
\]
where the kernels $\beta:\R^d\times\R^d\to\R^d$, $\varsigma:\R^d\times\R^d\to[0,\infty)$ and $\kappa:\R^d\times\R^d\to\R^{d\times d}$ are all bounded by $M$ and $\ell$-Lipschitz on $\R^d\times\R^d$. Then
\[
\Biggl( \E \Biggl[ \sup_{t\le T} |X^{i,N}_t-\bar X^i_t|^p \Biggr] \Biggr)^{1/p} \le \frac{C}{\sqrt{N}},
\]
with $C$ depending only on $d,s,p,T,M,\ell$, and not on $\mu_0$.
\end{theorem}

\begin{proof}
By the bound
\[
\Wr \Biggl( \frac{1}{N}\sum_j\delta_{x^j}, \frac{1}{N}\sum_j\delta_{y^j} \Biggr) \le \frac{1}{N}\sum_j|x^j-y^j|,
\]
the particle system has Lipschitz coefficients, so it is well posed by \cite[Thm.~V.7]{Protter2005}, and it is exchangeable. Let $D^i := X^{i,N}-\bar X^i$ and $V_r := \Wr(\bar\mu^N_r,\mu_r) \le 1$. As in the proof of Theorem~\ref{thm:chaos}, for almost every $r$ the differences of the coefficients are bounded by
\[
K \Biggl( |D^i_{r-}| + \frac{1}{N}\sum_j|D^j_{r-}| + V_r \Biggr).
\]

We work on the filtration generated by all the inputs $(X^i_0,W^i,L^i)_{i\ge1}$; the vector $(X^{i,N},\bar X^i)_{i\le N}$ is a permutation-equivariant functional of these i.i.d.\ inputs, hence jointly exchangeable. Let $T_0\in(0,1]$ be such that $T/T_0$ is an integer, and $t_k:=kT_0$. On $[t_k,t_{k+1}]$ we localise with the symmetric stopping time
\[
\tau_n := \inf \Biggl\{ t\ge t_k : \sum_j|D^j_t|\ge n \Biggr\} \wedge t_{k+1},
\]
and set $a_n := \E [ \sup_{t_k\le t\le\tau_n} |D^1_t|^p ]$. By convexity and exchangeability,
\[
\E \Biggl[ \sup_{t_k\le t\le\tau_n} \Biggl( \frac{1}{N}\sum_j|D^j_t| \Biggr)^p \Biggr] \le a_n.
\]
Lemma~\ref{lem:Lp}, together with
\begin{align*}
    \Biggl( \int_{t_k}^{t_{k+1}} V_r\,\dd r \Biggr)^p &\le T_0^{p-1} \int_0^T V_r^p\,\dd r, \\[2ex]
    \E \Biggl[ \Biggl( \int_{t_k}^{t_{k+1}} V_r^2\,\dd r \Biggr)^{p/2} \Biggr] &\le \Biggl( \E \Biggl[ \int_0^T V_r^2\,\dd r \Biggr] \Biggr)^{p/2},
\end{align*}
then gives, for $T_0$ small enough depending only on $d,s,p,K$ (by a condition analogous to \eqref{eq:T0}; note that $a_n<\infty$, since $\E|D^1_{t_k}|^p<\infty$ by the previous step and the localised integrands are bounded),
\[
a_n \le C_* \, \E \bigl[ |D^1_{t_k}|^p \bigr] + \frac{1}{2}a_n + C\,G_N,
\qquad
G_N := \E \Biggl[ \int_0^T V_r^p\,\dd r \Biggr] + \Biggl( \E \Biggl[ \int_0^T V_r^2\,\dd r \Biggr] \Biggr)^{p/2},
\]
with $C$ depending only on $d,s,p,K$. Letting $n\to\infty$,
\[
\E \Biggl[ \sup_{[t_k,t_{k+1}]} |D^1|^p \Biggr] \le 2C_* \, \E \bigl[ |D^1_{t_k}|^p \bigr] + 2C\,G_N,
\]
and since $D^1_0=0$, induction over the $\lceil T/T_0\rceil$ intervals gives $\E [ \sup_{[0,T]} |D^1|^p ] \le C' G_N$. The number of intervals depends only on $T$ and $T_0$, so $C'$ depends only on $d,s,p,K,T$, and not on $N$ or $\mu_0$. As in the proof of Theorem~\ref{thm:chaos}, $V_r\to0$ almost surely for every $r$ and $V\le1$, so $G_N\to0$ by dominated convergence.

In the kernel case, each coefficient difference splits as in the proof of Theorem~\ref{thm:chaos-rate} into a term bounded by
\[
\ell \Biggl( |D^i_{r-}| + \frac{1}{N}\sum_j|D^j_{r-}| \Biggr)
\]
and a fluctuation $R^i_r$ with $\E [ |R^i_r|^2 ] \le 4M^2/N$. The argument above applies with $V_r$ replaced by $|R^i_r|$, and since
\[
\E \bigl[ |R^i_r|^p \bigr] \le \Bigl( \E \bigl[ |R^i_r|^2 \bigr] \Bigr)^{p/2} \le \Biggl( \frac{4M^2}{N} \Biggr)^{p/2},
\]
we get
\[
\E \Biggl[ \sup_{[0,T]} |D^1|^p \Biggr] \le \frac{C}{N^{p/2}}.
\]
\end{proof}

The following example shows that Theorem~\ref{thm:jump-wp} fails for every $s\in(0,\frac12)$, even for a bounded coefficient that is Lipschitz with respect to $\Wr$ and depends only on the law. Cavallazzi \cite[Rem.~2]{Cavallazzi2025} shows that uniqueness may fail for coefficients that are Lipschitz in $W_\beta$ with $\beta<1$, through the deterministic equation $\dot y=|y|^\beta$ obtained with the unbounded drift $\displaystyle \int_{\R^d} |x|^\beta\,\mu(\dd x)$,
and no noise; see also \cite[Rem.~1.3]{DengHuang2025}. In Proposition~\ref{prop:counter}, the coefficient is bounded and Lipschitz for the bounded distance $\Wr$, and the non-uniqueness is produced by the heavy tails of the stable noise itself.

\begin{proposition}\label{prop:counter}
Let $s\in(0,\frac{1}{2})$, $\alpha := 2s$, and
\[
\gamma(\mu) := \Biggl( \int_{\R^d} \bigl( |x|\wedge1 \bigr) \, \mu(\dd x) \Biggr) I_d \in \R^{d\times d}.
\]
Then $\gamma$ satisfies Assumption~\ref{ass:Hg} with $b=\Sigma=0$ and $K=\sqrt{d}$, since
\[
|\gamma(\mu)-\gamma(\nu)| \le \sqrt{d} \, \Wr(\mu,\nu),
\]
and the equation
\[
\dd X_t = \gamma(\mu_t)\,\dd L_t, \qquad \mu_t = \Law(X_t), \qquad X_0 = 0,
\]
has two different solutions: the trivial solution $X\equiv0$, and a solution with
\[
\Law(X_t) = \Law \Bigl( S_t^{1/\alpha} L_1 \Bigr),
\]
where $S$ is a continuous function with $S_0=0$ and $S_t>0$ for $t>0$.
\end{proposition}

\begin{proof}
Since $x\mapsto|x|\wedge1$ is bounded by $1$ and $1$-Lipschitz, we have $\big||x|\wedge1-|y|\wedge1\big| \le |x-y|\wedge1$, and the Lipschitz bound follows by integrating against an optimal coupling. The process $X\equiv0$ is a solution, since $\gamma(\delta_0)=0$.

For $v\ge0$ we set $m(v) := \E\bigl[|v^{1/\alpha}L_1|\wedge1\bigr]$ and $g(v) := m(v)^\alpha$. By dominated convergence, $g$ is continuous; moreover $g(0)=0$, $g(v)>0$ for $v>0$, and $g\le1$. The first coordinate of $L_1$ is a symmetric $\alpha$-stable random variable, so $\Prob(|L_1|>u) \ge c\,u^{-\alpha}$ for $u\ge u_0$, with $c>0$ \cite[Property~1.2.15]{SamorodnitskyTaqqu1994}. Hence, for $v$ small,
\[
m(v) \ge \Prob \bigl( |L_1|>v^{-1/\alpha} \bigr) \ge cv, \qquad \text{so that} \quad g(v) \ge c^\alpha v^\alpha.
\]
Since $\alpha<1$, the function
\[
G(v) := \int_0^v \frac{\dd u}{g(u)}, \qquad v\ge0,
\]
is finite, continuous and strictly increasing, with $G(0)=0$ and $G(v)\to\infty$ as $v\to\infty$ because $g\le1$. Its inverse $S_t := G^{-1}(t)$, $t\ge0$, satisfies $S_0=0$, $S_t>0$ for $t>0$, and $\displaystyle S_t = \int_0^t g(S_r)\,\dd r$. Note that $S\equiv0$ also solves this last equation; the two solutions $0$ and $S$ give the two solutions of the McKean--Vlasov equation.

We now put $h_r := m(S_r)$, a deterministic continuous function, and $\displaystyle X_t := \int_0^t h_r\,\dd L_r$.
Approximating $h$ by step functions and using the independence and stationarity of the increments of $L$, we get
\[
\E \Bigl[ e^{\iu\xi\cdot X_t} \Bigr] = \exp \Biggl( -|\xi|^\alpha \int_0^t h_r^\alpha\,\dd r \Biggr) = \exp \bigl( -|\xi|^\alpha S_t \bigr),
\]
since $h_r^\alpha=g(S_r)$. Hence $\Law(X_t)=\Law(S_t^{1/\alpha}L_1)$, and $\gamma(\Law(X_t)) = m(S_t)I_d = h_tI_d$. Therefore $X$ solves $\dd X_t=\gamma(\mu_t)\,\dd L_t$ with $X_0=0$, and $\Law(X_t)\neq\delta_0$ for $t>0$. Both solutions are defined on the same probability space and driven by the same process $L$, so pathwise uniqueness and uniqueness in law both fail.
\end{proof}

\begin{remark}
The proof rests on a classical fact: the equation $\dot S=g(S)$, $S_0=0$, has a solution other than $S\equiv0$ if and only if $\displaystyle\int_0^{1/2}\frac{\dd u}{g(u)}<\infty$, that is, if $g$ is large enough near $0$. For $s<\frac12$ this holds because $g(v)\ge c^\alpha v^{2s}$ with $2s<1$. For $s\ge\frac12$ it fails, because $g$ is too small near $0$:
\begin{itemize}
\item if $s>\frac12$, then $\E|L_1|<\infty$ and $m(v)\le v^{1/\alpha}\,\E|L_1|$, hence $g(v)\le (\E|L_1|)^\alpha\, v$;
\item if $s=\frac12$, the tail bound $\Prob(|L_1|>u)\le C/u$ gives $m(v)\le C' v\log(1/v)$ for small $v$, and here $g=m$.
\end{itemize}
In both cases the integral diverges, so $S\equiv0$ is the only solution and the construction gives nothing. For $s>\frac12$ this agrees with Theorem~\ref{thm:jump-wp}. For $s=\frac12$, the question remains open.
\end{remark}

\appendix

\section{Linear dynamics: generator, symbol and law}
\label{sec:linear}

This appendix collects the linear facts used in Sections~\ref{sec:nonlinear} and~\ref{sec:jumps},
where they are applied to the dynamics with frozen law. They are classical;
we give short proofs to fix the conventions, in particular the signs.

\begin{assumption}\label{ass:lin}
$b:\R^d\to\R^d$ and $\sigma:\R^d\to[0,\infty)$ are such that $b$ and
$\Sigma:=\sqrt{2\sigma}$ are bounded and Lipschitz continuous.
\end{assumption}

Let $W$ be a standard $d$-dimensional $(\cF_t)$-Brownian motion and $L$ an
independent isotropic $2s$-stable $(\cF_t)$-L\'evy process. Consider
\begin{equation}\label{eq:SDE-lin}
  \dd X_t=b(X_t)\,\dd t+\Sigma(X_t)\,\dd W_t+\dd L_t,\qquad X_0\ \
  \cF_0\text{-measurable},
\end{equation}
and the \emph{mixed generator}
\begin{equation}\label{eq:gen-mixed}
  \cA f(x):=b(x)\cdot\nabla f(x)+\sigma(x)\,\Delta f(x)-(-\Delta)^s f(x),
  \qquad f\in C_b^2(\R^d).
\end{equation}
Its \emph{carr\'e du champ}
$\Gamma(f,g):=\cA(fg)-f\cA g-g\cA f$ \cite{BGL2014} is computed directly:
the drift terms cancel because $b\cdot\nabla$ is a derivation, and the
compensator terms in \eqref{eq:fraclap} cancel for the same reason. This gives
\begin{equation}\label{eq:carre}
  \Gamma(f,g)(x)=2\sigma(x)\,\nabla f(x)\cdot\nabla g(x)
  +\int_{\R^d\setminus\{0\}}\big(f(x+y)-f(x)\big)\big(g(x+y)-g(x)\big)\,
  \nu_s(\dd y).
\end{equation}

Since the jump noise in \eqref{eq:SDE-lin} is additive, \eqref{eq:SDE-lin}
can be written as $\displaystyle X_t=J_t+\int_0^tF(X)_{r-}\,\dd Z_r$ with $J:=X_0+L$,
$Z:=(t,W_t)$ and $F(X)_r:=(b(X_r),\Sigma(X_r))$. The operator $F$ is
functional Lipschitz. By \cite[Thm.~V.7]{Protter2005}, \eqref{eq:SDE-lin}
has a unique c\`adl\`ag adapted solution for every $\cF_0$-measurable $X_0$,
without any integrability requirement on $X_0$. We write $\E^x$ when
$X_0=x$.

\begin{proposition}[Pathwise reading: It\^o's formula]\label{prop:ito}
Under Assumption~\ref{ass:lin}, for every $\varphi\in C_b^2(\R^d)$,
\begin{equation}\label{eq:ito-lin}
  \varphi(X_t)=\varphi(X_0)+\int_0^t\cA\varphi(X_r)\,\dd r+M^\varphi_t,
\end{equation}
where
\[
  M^\varphi_t:=\int_0^t\Sigma(X_r)\nabla\varphi(X_r)\cdot\dd W_r
  +\int_0^t\!\!\int_{\R^d\setminus\{0\}}\big[\varphi(X_{r-}+y)-\varphi(X_{r-})\big]
  \tilde N(\dd r,\dd y)
\]
is a square-integrable martingale with predictable quadratic variation
$\displaystyle \langle M^\varphi\rangle_t=\int_0^t\Gamma(\varphi,\varphi)(X_r)\,\dd r$.
\end{proposition}

\begin{proof}
The continuous martingale part of $X$ is $\displaystyle \int\Sigma(X_r)\,\dd W_r$, with
$\dd[X^{c,i},X^{c,j}]_r=\Sigma(X_r)^2\delta_{ij}\,\dd r=2\sigma(X_r)\delta_{ij}\,\dd r$.
It\^o's formula for semimartingales with jumps
\cite[Thm.~4.4.7]{Applebaum2009}, combined with \eqref{eq:LI}, gives
\begin{align*}
  \varphi(X_t)-\varphi(X_0)
  &=\int_0^t b\cdot\nabla\varphi(X_r)\,\dd r
   +\int_0^t\Sigma\nabla\varphi(X_r)\cdot\dd W_r
   +\int_0^t\sigma\Delta\varphi(X_r)\,\dd r\\
  &\quad+\int_0^t\!\!\int_{\{|y|<1\}}\big[\varphi(X_{r-}+y)-\varphi(X_{r-})\big]\tilde N(\dd r,\dd y)\\
  &\quad+\int_0^t\!\!\int_{\{|y|\ge1\}}\big[\varphi(X_{r-}+y)-\varphi(X_{r-})\big]N(\dd r,\dd y)\\
  &\quad+\int_0^t\!\!\int_{\{|y|<1\}}\big[\varphi(X_{r}+y)-\varphi(X_{r})-\nabla\varphi(X_r)\cdot y\big]\nu_s(\dd y)\,\dd r .
\end{align*}
Note that the jump term carries the correction $-\nabla\varphi\cdot\Delta X$,
so no finite-variation assumption on the jumps is needed (this matters when
$2s\ge1$). Since $\nu_s(\{|y|\ge1\})<\infty$ and the integrand is bounded,
we may write $N=\tilde N+\dd r\,\nu_s(\dd y)$ on $\{|y|\ge1\}$. Collecting
the $\dd r$-terms and using \eqref{eq:fraclap} yields \eqref{eq:ito-lin}.
We used that $X_{r-}=X_r$ for all but countably many $r$.

The Brownian integrand is bounded, since $\Sigma$ and $\nabla\varphi$ are.
The Poisson integrand satisfies
$\displaystyle \int|\varphi(x+y)-\varphi(x)|^2\nu_s(\dd y)\le\int\min(\norm{\nabla\varphi}_\infty^2|y|^2,4\norm\varphi_\infty^2)\,\nu_s(\dd y)<\infty$,
uniformly in $x$. Hence both stochastic integrals are square-integrable
martingales, orthogonal to each other, with predictable quadratic variations
$\displaystyle \int2\sigma|\nabla\varphi|^2(X_r)\,\dd r$ and
$\displaystyle \int\!\int(\varphi(X_r+y)-\varphi(X_r))^2\nu_s(\dd y)\,\dd r$. Their sum is
$\displaystyle \int\Gamma(\varphi,\varphi)(X_r)\,\dd r$ by \eqref{eq:carre}.
\end{proof}

\begin{proposition}[Analytic reading: symbol]\label{prop:analytic}
Under Assumption~\ref{ass:lin}, the symbol \eqref{eq:symbol} of $\cA$ is
\begin{equation}\label{eq:symbol-mixed}
  p(x,\xi)=-\iu\, b(x)\cdot\xi+\sigma(x)|\xi|^2+|\xi|^{2s},
\end{equation}
and for every $f\in\Sch(\R^d)$, $\cA f=-p(x,D)f$ in the sense of
\eqref{eq:pdo}. If $b$ and $\sigma$ are constant, $p(x,\xi)=\psi(\xi)$ is
the characteristic exponent of the L\'evy process
$bt+\sqrt{2\sigma}W_t+L_t$, whose triplet is $(b,2\sigma I_d,\nu_s)$.
\end{proposition}

\begin{proof}
We have $\nabla e_\xi=\iu\xi e_\xi$ and $\Delta e_\xi=-|\xi|^2e_\xi$, and
\eqref{eq:fraclap-plane} handles the jump part. Hence
$\cA e_\xi(x)=\big(\iu\, b(x)\cdot\xi-\sigma(x)|\xi|^2-|\xi|^{2s}\big)e_\xi(x)$,
which gives \eqref{eq:symbol-mixed}. For $f\in\Sch$, write
$\displaystyle f=(2\pi)^{-d}\int\hat f(\xi)e_\xi\,\dd\xi$. The map
$\xi\mapsto\hat f(\xi)e_\xi$ is Bochner integrable in $C_b^2$, since
$\norm{e_\xi}_{C_b^2}\le C(1+|\xi|)^2$ and $\hat f\in\Sch$. Moreover, $\cA$
is a bounded operator $C_b^2\to C_b$ (Definition~\ref{def:fraclap}), so it
commutes with the integral. The last statement follows from
Example~\ref{ex:basic}, Definition~\ref{def:stable} and independence.
\end{proof}

\begin{proposition}[Probabilistic formulas for $\cA$ and $p$]\label{prop:schilling}
Under Assumption~\ref{ass:lin}, for every $\varphi\in C_b^2(\R^d)$, every
$x\in\R^d$ and every $\xi\in\R^d$,
\begin{equation}\label{eq:schilling}
  \cA\varphi(x)=\lim_{t\downarrow0}\frac{\E^x\varphi(X_t)-\varphi(x)}{t},
  \qquad
  p(x,\xi)=-\lim_{t\downarrow0}\frac{\E^x\big[e^{\iu\xi\cdot(X_t-x)}\big]-1}{t}.
\end{equation}
\end{proposition}

\begin{proof}
Taking expectations in \eqref{eq:ito-lin} gives
$\displaystyle t^{-1}(\E^x\varphi(X_t)-\varphi(x))=t^{-1}\int_0^t\E^x[\cA\varphi(X_r)]\,\dd r$.
The function $\cA\varphi$ is bounded and continuous, and $X_r\to x$ as
$r\downarrow0$ by right-continuity. By dominated convergence,
$r\mapsto\E^x[\cA\varphi(X_r)]$ is right-continuous at $0$, which gives the
first formula. For the second, apply the first formula to the real and
imaginary parts of $e_\xi$, use $\cA e_\xi=-p(\cdot,\xi)e_\xi$, and multiply
by $e_{-\xi}(x)$.
\end{proof}

The second formula in \eqref{eq:schilling} is the probabilistic definition
of the symbol of a Feller process due to Schilling; see
\cite{SchillingSchnurr2010} for SDEs driven by L\'evy processes.

Taking expectations in \eqref{eq:ito-lin} shows that $\mu_t:=\Law(X_t)$
solves the forward Kolmogorov equation
$\displaystyle \ip{\mu_t}{\varphi}=\ip{\mu_0}{\varphi}+\int_0^t\ip{\mu_r}{\cA\varphi}\,\dd r$
for $\varphi\in C_b^2(\R^d)$, and testing it against $e_\xi$ with
$\mu_0=\delta_x$ recovers the symbol through \eqref{eq:schilling}. In this
way, the generator is described by its action on paths
(Proposition~\ref{prop:ito}), on plane waves
(Proposition~\ref{prop:analytic}) and on laws.
Section~\ref{sec:nonlinear} extends each of these descriptions to
McKean--Vlasov flights.

\end{document}